\documentclass[12pt]{article}
\usepackage{makeidx}
\usepackage{times}
\usepackage{amssymb}
\usepackage{amsmath}
\numberwithin{equation}{section}

\usepackage{color}
\usepackage{graphicx}
\usepackage{rotating}
\usepackage{pdflscape}
\usepackage{epstopdf}
\usepackage[round,authoryear,comma]{natbib}
\usepackage{natbib,hyperref}
\usepackage{booktabs}
\usepackage{setspace}
\usepackage{pgfplots}
\usepackage{subcaption}
\usepackage{xcolor}
\usepackage{multirow}
\hypersetup{
	colorlinks   = true, 
	urlcolor     = blue, 
	linkcolor    = blue, 
	citecolor   = blue 
}

\providecommand{\U}[1]{\protect\rule{.1in}{.1in}}
\newtheorem{theorem}{Theorem}[section]

\newtheorem{definition}{Definition}[section]

\newtheorem{lemma}{Lemma}[section]

\newtheorem{proposition}{Proposition}
\newtheorem{remark}[theorem]{Remark}

\parskip=\medskipamount

\begin{document}

\title{Monotone tail functions: definitions, properties, and application to risk-reducing strategies}
\date{  }
\author{Hamza Hanbali \\ \small Department of Econometrics and Business Statistics, Monash University, Australia \\  \texttt{{\small hamza.hanbali@monash.edu}} \bigskip
	\and Dani\"{e}%
l Linders \\ \small
Department of Mathematics, University of Illinois at Urbana-Champaign, USA \\ \small Faculty of Economics and Business, University of Amsterdam, Netherlands. \\ \texttt{\small dlinders@illinois.edu}}
\maketitle
\begin{abstract}
This paper studies properties of functions having monotone tails. We extend Theorem 1 of \cite{Dhaene2002a} and show how the tail quantiles of a random variable transformed with a monotone tail function can be expressed as the transformed tail quantiles of the original random variable. The main result is intuitive, in that \cite{Dhaene2002a}'s properties still hold, but only for certain quantile values. However, the proof presents some complications that arise especially when the function involved has discontinuities. 

We consider different situations where monotone tail functions occur and can be useful, such as the evaluation of the payoff of option trading strategies and the present value of insurance contracts providing both death and survival benefits. The paper also applies monotone tail functions to study quadrant perfect dependence, and shows how this dependence structure integrates within the framework of monotone tail functions. Moreover, we apply the theory to the problem of risk reduction and investigate conditions on a hedger ensuring efficient reductions of required economic capital.

\end{abstract}
\newpage 
\section{Introduction}
\label{Introduction} 

A standard result derived in \cite{Dhaene2002a} expresses the quantile of a random variable transformed by a monotone function as the transform of the original random variable's quantile. The present paper derives sufficient conditions such that this result remains to hold when the transformation involves a function which is not monotone. We extend Theorem 1 of \cite{Dhaene2002a} by investigating \textit{monotone tail} functions $h$: functions which are monotone after or before a certain threshold. 
We show in this paper that if $h$ has a monotone tail, then transforming the original random variable using $h$ and then taking the quantile of the transformed random variable, is equal to first taking the quantile and then transforming this quantile using the function $h$, provided the confidence level $p$ of the quantile is within an explicitly given interval. This result complies with the intuition, but the proof presents some challenges when the function $h$ has discontinuities.

The paper studies different situations in which such functions having a monotone tail arise. Option trading strategies involving combinations of calls and puts are typical examples where the payoff function of the strategy has a monotone tail. These strategies are often used to neutralize portfolios in some parts of the payoff's support. Determining the quantiles in the tail of the payoff's distribution can be achieved using the theory developed in this paper. We also show that monotone tail functions arise in the distribution of the payout of insurance contracts such as endowment assurances. These contracts pay a benefit in case of death or survival of the policyholder. Thus, the present value of the benefit is expressed as the sum of a non-decreasing function and a non-increasing function of the remaining lifetime, and the resulting function has a monotone tail.

A related concept is \textit{quadrant perfect dependence}. For example, upper comonotonicity, which was introduced in \cite{Cheung2009}, refers to a perfect positive dependence structure in the upper quadrant of their joint support; see also \cite{DongCheungYang}, \cite{ChungLo} and \cite{Cheungetal}. In general, quadrant perfect dependence refers to a dependence structure under which two random variables are perfectly positive or negative dependent in a quadrant. The possible generalization of upper comonotonicity to other quadrants has already been discussed in e.g. \cite{DongCheungYang}. The present paper studies these dependence structures through the lens of the proposed monotone tail functions. We show that monotonicity in the tail of the difference of the two quantile functions of two comonotonic random variables $R_1$ and $R_2$ is equivalent with quadrant perfect dependence of the random vector $(R_2, R_1-R_2)$. 

The paper also employs monotone tail functions and quadrant perfect dependence to address the problem of selecting an appropriate Value-at-Risk reducer. Risk reduction consists in finding an appropriate asset $R_2$ whose interaction with a given liability $R_1$ is such that the risk of the residual liability $R_1-R_2$ is reduced; see \cite{Ahnetal} and \cite{DeelstraEzzineHeymanVanmaele} for some early studies. In terms of solvency capital, an insurance company faces the problem of selecting a hedger $R_2$ such that the amount of required regulatory capital is reduced, where this capital is measured using an appropriate quantile, i.e.\ Value-at-Risk (VaR). Intuitively, a candidate risk reducer $R_2$ should be (strongly) positively dependent with $R_1$. But even when $R_1$ and $R_2$ are perfectly positively dependent (i.e. comonotonic), selecting an adequate asset to reduce the risk of the liabilities is not a straightforward task. We work under a setting where $R_1$ and $R_2$ are comonotonic random variables, i.e.\ both random variables are non-decreasing functions of the same random source. This assumption provides a useful proxy when studying risk reducers; see also \cite{Ahnetal} and \cite{CheungDhaeneLoTang}. In this case, the difference $Z=R_1-R_2$ can be expressed as $Z\overset{d}{=}h(U)$, where $U$ is a random variable uniformly distributed over the unit interval and 
\begin{equation}\label{intro-1}
h(p)=F^{-1}_{R_1}(p)-F^{-1}_{R_2}(p).
\end{equation}
This equality does not necessarily mean that the VaR of the residual liability $Z$ is the difference between the VaR of the initial position and that of the candidate VaR-reducer. Theorem 1 of \cite{Dhaene2002a} ensures that such a decomposition holds when the function $h$ is monotone. In our case, however, the function $h$ is not necessarily monotone. As a consequence, we cannot directly express the quantile of the difference as a difference of quantiles. When using the VaR to determine solvency capital, one is mainly interested in quantiles with sufficiently high level $p$. The theory on functions having a monotone tail is illustrated by showing that the decomposition of the quantile of $Z$ holds if we limit the range of $p$ to some interval of interest, especially in the tail. Alternatively, it should be acknowledged that the decomposition of the quantile of the difference can be derived using the theory on upper comonotonicity and the extension of quadrant perfect dependence. Either way, this decomposition allows to derive conditions on the choice of an appropriate VaR-reducer. This leads to the conclusion that if the upper tail of the residual liability is dominated by the initial liability, the intended reduction in solvency capital is achieved provided the asset $R_2$ is traded at a reasonable price. In the opposite case, a reasonable price for $R_2$ might not exist.

In order to gain more insights into the comonotonic risk mitigation process, we consider a general setting where $R_1$ and $R_2$ have equally distributed standardized log-returns. We show in this case that the support of the residual liability always has a monotone upper tail. We also find that the price conditions for VaR reduction depend on the relative standard deviations of the candidate VaR reducer and the initial liability. On the one hand, in case the standard deviation associated with $R_2$ is smaller than that of the initial liability, we find that $R_2$ is a VaR-reducer provided its price is smaller than its VaR at the confidence level $p$. Since in practice $p$ takes high values, this condition will be satisfied for any reasonable price of $R_2$. On the other hand, in case the standard deviation associated with the candidate risk reducer is larger than that of the initial liability, we find that the losses from the residual liability are bounded, and it can still occur that $R_2$ is not a VaR-reducer. What explains this result is that in such situations, the price of the candidate hedger is too large compared to the capital relief, and the capital required to achieve the risk reduction exceeds the capital needed to cover the liability $R_1$ when it is not hedged at all.
 
Comonotonicity of the initial liability and the risk reducer implies that the risk reduction is efficient, in the sense that a large loss is followed by a large payoff of the risk reducer. However, if the payoff from the risk reducer is too large, the initial claim is over-hedged, which can potentially be costly and therefore not necessarily desirable. One situation which avoids this over-hedging is when the support of the residual loss has a non-decreasing upper tail. This ensures that if the risk reducer is large, the residual loss will also increase. In other words, the risk reducer does not dominate the initial liability in the upper tail. As we show in this paper, this is equivalent with the situation where the risk reducer and the residual liabilities are upper comonotonic.

An alternative, and more general, characterization of risk mitigation was given in \cite{CheungDhaeneLoTang}, who define the notion of a risk reducer in terms of the Tail Value-at-Risk (TVaR). They show that if $R_1$ and $R_2$ move in the same direction, and hence $R_2$ would presumably compensate the high losses of $R_1$, the hedged portfolio is not necessarily less risky than the initial position. In other words, there remain situations where not hedging at all is better than hedging with a comonotonic asset. Nevertheless, when the asset and the residual liability are comonotonic, then this asset is a risk reducer in the sense of \cite{CheungDhaeneLoTang}. We also refer to \cite{HeTangZhang_2016} for a generalization to risk reducers which are not perfectly dependent. Our analysis shows that the conclusions on VaR-reducers hold for TVaR-reducers too, if we limit the condition of TVaR reduction to some values of the confidence level.

The remainder of the paper is organized as follows. Section \ref{Section-Preliminaries} contains preliminary definitions and notations. Section \ref{Section-Quantiles} provides the main technical results of the paper, with the definitions of all four types of monotone tail functions and the theorem extending that of \cite{Dhaene2002a} in each case. We study monotone tail functions in the context of option trading strategies and endowment assurances in Section \ref{Section-Results}.  Four types of quadrant perfect dependence are discussed and their relationship with monotone tail functions is established in Section \ref{SectionQuadrantDependence}.  We address in Section \ref{Section-CapitalReducers} the problem of risk reduction using monotone tail functions, and we investigate conditions on the distribution of $R_2$ leading to a VaR reduction. We conclude in Section \ref{Section-Final}. All proofs can be found in Appendix \ref{Section-Appendix}.


\section{Preliminaries} \label{Section-Preliminaries}
All random variables are defined on a common probability space $\left(\Omega,\mathcal{F},
\mathbb{P} \right)$ and can be discrete or continuous. The cumulative distribution function (cdf) of a random variable
$X$ is denoted by $F_X$, where $X$ takes values in the real interval $\mathcal{I}\subseteq \mathbb{R}$. The left inverse cdf $F_X^{-1}$ of $X$ is defined as
\begin{equation}
F_{X}^{-1}\left(  p\right)  =\inf\left\{  x\in\mathbb{R}\mid F_{X}(x)\geq
p\right\}  ,\qquad p\in\left[  0,1\right]  , \label{inverse}%
\end{equation}
with $\inf\emptyset=+\infty$, by convention, and the Value-at-Risk of $X$ at the level $p$, denoted by $\text{VaR}_p[X]$, corresponds to $F^{-1}_{X}(p)$. The Tail Value-at-Risk at the level $p$ is defined as follow:
$$
\text{TVaR}_p[X] = \frac{1}{1-p}\int_p^1 F^{-1}_X(q)\text{d}q.
$$
The Left TVaR at the level $p$ is defined as follows:
$$\text{LTVaR}_p[X] = \frac{1}{p}\int_0^p F^{-1+}_X(q)\text{d}q,$$
where the right inverse $F^{-1+}_X$ is defined as follows: 
\begin{equation}
F_{X}^{-1+}\left(  p\right)  =\sup\left\{  x\in\mathbb{R}\mid F_{X}(x)\leq
p\right\}  ,\qquad p\in\left[  0,1\right]  , \label{+inverse}%
\end{equation}
where $\sup\emptyset=-\infty$, by convention. We also use the generalized $\alpha$-inverse $F_{X}^{-1\left(  \alpha\right)  }$:
\begin{equation}
F_{X}^{-1\left(  \alpha\right)  }\left(  p\right)  =(1-\alpha) F_{X}^{-1}\left(
p\right)  +\alpha  F_{X}^{-1+}\left( p\right)  ,\qquad
p\in\left[ 0,1\right] \text{ and } \alpha \in[0,1], \label{alpha-inverse}%
\end{equation}
More details on generalized inverses can be found in \cite{Embrechts2013}.

A set $\mathcal{S}$ is comonotonic if for any couples $(x_1,x_2), (y_1,y_2)\in \mathcal{S}$, we have that $x_1\leq y_1$ and $x_2\leq y_2$, or $x_1\geq y_1$ and $x_2\geq y_2$. More specifically, a comonotonic set is non-decreasing in both its components. A random vector $\left(X_1,X_2\right)$ is said to be comonotonic if its support is a comonotonic set. Analogously, a random vector is said to be counter-monotonic if its support is a counter-monotonic set, where a set $\mathcal{S}$ is counter-monotonic if for any $(x_1,x_2), (y_1,y_2)\in \mathcal{S}$, we have that $x_1\leq y_1$ and $x_2\geq y_2$, or $x_1\geq y_1$ and $x_2\leq y_2$.

\section{Quantiles of random variables transformed with monotone tail functions} \label{Section-Quantiles}
\subsection{Monotone functions}
Consider a function $h: \mathcal{I}\mapsto \mathcal{J}$, with $\mathcal{I},\mathcal{J}\subseteq \mathbb{R}$. \cite{Dhaene2002a} derive conditions such that the quantiles of a transformed random variable coincide with the  transformed quantiles of the original random variable. This is formally stated in the following theorem. 
\begin{theorem} \label{Theorem-Dhaene}
	Consider the random variable $X$ with cdf $F_X$, and a function $h:\mathcal{I}\mapsto \mathcal{J}$.
	\begin{enumerate}
		\item If $h$ is non-decreasing and left-continuous, then: 
		\begin{equation}\label{T-h-D1}
			F_{h(X)}^{-1}(p)=h\left(F_X^{-1}(p) \right), \text{\qquad for} \ p\in (0,1).
		\end{equation}
		\item If $h$ is non-decreasing and right-continuous, then:
		\begin{equation}\label{T-h-D2}
			F_{h(X)}^{-1+}(p)=h\left(F_X^{-1+}(p) \right), \text{\qquad for} \ p\in (0,1).
		\end{equation}
		\item If $h$ is non-increasing and left-continuous, then: 
		\begin{equation}\label{T-h-D3}
			F_{h(X)}^{-1+}(p)=h\left(F_X^{-1}(1-p) \right), \text{\qquad for} \ p\in (0,1).
		\end{equation}
		\item If $h$ is non-increasing and right-continuous, then:
		\begin{equation}\label{T-h-D4}
			F_{h(X)}^{-1}(p)=h\left(F_X^{-1+}(1-p) \right), \text{\qquad for} \ p\in (0,1).
		\end{equation}
	\end{enumerate}
\end{theorem}

In this section, it is shown that, for some values of $p$, the results of Theorem \ref{Theorem-Dhaene} hold even in case $h$ is not monotone everywhere. The functions under interest are monotone after or before a certain threshold. Such functions are said to have a monotone tail. It is proven that Theorem \ref{Section-Quantiles} remains to hold for these functions, provided the probabilities $p$ are either sufficiently large, or sufficiently small. 

\subsection{Monotone tail functions}
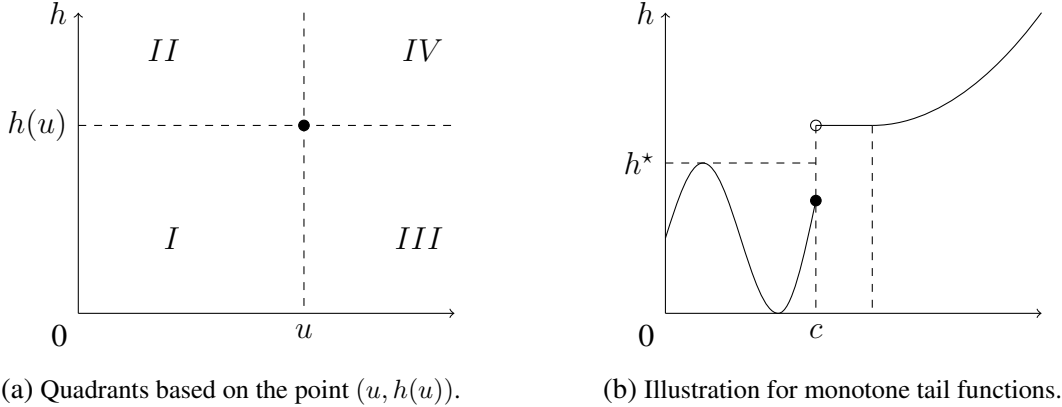
\begin{figure}[!h]
	\centering
	\begin{subfigure}{0.5\textwidth}
		\centering
		\begin{tikzpicture}
			\draw[->] (0,0) -- (5,0) node[anchor=north]{} ;
			\draw[->] (0,0) -- (0,4) node[anchor=east] {$h$};
			\draw	(0,0) node[anchor=north east] {0}
			(3,0) node[anchor=north] {$u$}
			(1.5,1) node[anchor=east] {$I$}
			(1.5,3.5) node[anchor=east] {$II$}
			(5,1) node[anchor=east] {$III$}
			(5,3.5) node[anchor=east] {$IV$}
			(0,2.5) node[anchor=east] {$h(u)$};
			\filldraw (3,2.5) circle (2pt);
			\draw[dashed] (0,2.5) -- (5,2.5); 
			\draw[dashed] (3,4) -- (3,0); 
		\end{tikzpicture}
		\caption{\footnotesize Quadrants based on the point $(u,h(u))$. }\label{Fig:Tails1a}
	\end{subfigure}%
	\begin{subfigure}{0.5\textwidth}
		\centering
		\begin{tikzpicture}
			\draw[->] (0,0) -- (5,0) node[anchor=north]{} ;
			\draw[->] (0,0) -- (0,4) node[anchor=east] {$h$};
			\draw	(0,0) node[anchor=north east] {0}
			(2,0) node[anchor=north] {$c$}
			(0,2) node[anchor=east] {$h^{\star}$};
			\draw (0,1) sin (0.5,2) cos (1,1) sin (1.5,0) cos (2,1.5);
			\draw (2,2.5) -- (2.75,2.5);
			\draw (2.75,2.5) parabola (5,4) ;
			\filldraw (2,1.5) circle (2pt);
			\draw (2,2.5) circle (2pt);
			\draw[dashed] (0,2) -- (2,2); 
			\draw[dashed] (2,2.5) -- (2,0); 
			\draw[dashed] (2.75,2.5) -- (2.75,0); 
		\end{tikzpicture}
		\caption{\footnotesize Illustration for monotone tail functions. }\label{Fig:Tails1}
		
	\end{subfigure}
	\caption{\footnotesize Functions with a non-decreasing upper tail.}
\end{figure}

Generally speaking, a function $h:\mathcal{I}\mapsto \mathcal{J}$ is said to have a monotone tail if $\mathcal{I}\times \mathcal{J}$ can be divided into four quadrants (see Figure \ref{Fig:Tails1a}), such that $h$ is monotone in one of the quadrants. Naturally, there are four types of monotone tail functions. The function $h$ can be either non-decreasing or non-increasing \textit{after} a certain threshold. Such functions are said to have a monotone upper, or right, tail. Analogously, the function $h$ can be either non-decreasing or non-increasing \textit{before} a certain threshold. Such functions are said to have a monotone lower, or left, tail.  

In the sequel, the function used to transform the random variable $X$ is assumed to be real-valued, and can only have right or left continuous jumps. The notation $h$ is used for functions with non-decreasing tail, and $g$ for functions with non-increasing tail. Further, the subscripts ``$u$'' and ``$l$'' are used to indicate whether the monotonicity is in the upper or the lower tail, respectively. The initial definition is that of functions having a non-decreasing upper tail.
\begin{definition}[non-decreasing upper tail]\label{Def:Tailmonotonic}
	The function $h_u$ is said to have a \emph{non-decreasing upper tail} if there exists a constant $x'\in \mathcal{I}$ such that: 
	\begin{enumerate}
		\item $h_u(x)\leq h_u(x')$, for all elements of the set $\left\{ x\in \mathcal{I}|\ x<x'\right\}$, provided it is not empty,
		\item $h_u(x_1)\leq h_u(x_2)$, for all $x_1, x_2 \in \mathcal{I}$ with $x' \leq x_1 \leq x_2$.
	\end{enumerate}
	The function $h_u$ is said to have a \emph{non-increasing upper tail with threshold} $c_u=\inf \mathcal{H}_u$, where $\mathcal{H}_u$ contains all $x'$ satisfying Conditions (1) and (2).
\end{definition}

Condition (1) of Definition \ref{Def:Tailmonotonic} states that the function $h_u$ only takes values in quadrants I and IV of Figure \ref{Fig:Tails1a}, and does not enter the quadrants II and III. Condition (2) of Definition \ref{Def:Tailmonotonic} states that the function should be non-decreasing in quadrant IV.

Note that when $h_u$ is non-decreasing, then $\mathcal{H}_u=\mathcal{I}$. In this case, it is sufficient to use the result of \cite{Dhaene2002a} stated in Theorem \ref{Theorem-Dhaene}. When the function $h_u$ does not have a non-decreasing upper tail, the set $\mathcal{H}_u$ is empty, and by convention, $c_u=\sup \mathcal{I}$. In this case, monotone tail functions are not relevant. Thus, throughout this paper, these trivial cases are ruled out, and two assumptions are used, namely, $\mathcal{H}_u$ is not empty, and $c_u\in \left(\inf \mathcal{I},\sup \mathcal{I}\right)$.

The other forms of monotone tail functions are now defined.
\begin{definition}[non-increasing upper tail]\label{Def:Tailmonotonic_nonincreasing}
	The function $g_u$ is said to have a \emph{non-increasing upper tail} if there exists $x'\in \mathcal{I}$ such that: 
	\begin{enumerate}
		\item $g_u(x)\geq g_u(x')$, for all elements of the set $\left\{ x\in \mathcal{I}|\ x<x'\right\}$, provided it is not empty,
		\item $g_u(x_1)\geq g_u(x_2)$, for all $x_1, x_2 \in \mathcal{I}$ with $x'\leq x_1 \leq x_2$.
	\end{enumerate}
	The function $g_u$ is said to have a \emph{non-increasing upper tail with threshold} $d_u=\inf \mathcal{G}_u$, where $\mathcal{G}_u$ contains all $x'$ satisfying Conditions (1) and (2).
\end{definition}

\begin{definition}[non-decreasing lower tail] \label{Def:LowerTailmonotonic}
	The function $h_l$ is said to have a \emph{non-decreasing lower tail} if there exists $x'\in \mathcal{I}$ such that: 
	\begin{enumerate}
		\item $h_l(x)\geq h_l(x')$, for all elements of the set $\left\{ x\in \mathcal{I}|\ x>x'\right\}$, provided it is not empty,
		\item $h_l(x_1)\leq h_l(x_2)$, for all $x_1, x_2 \in \mathcal{I}$ with $x_1 \leq x_2\leq x'$.
	\end{enumerate}
	The function $h_l$ is said to have a \emph{non-decreasing lower tail with threshold} $c_l=\sup \mathcal{H}_l$, where $\mathcal{H}_l$ contains all $x'$ satisfying Conditions (1) and (2).
\end{definition}

\begin{definition}[non-increasing lower tail] \label{Def:LowerTailmonotonic_nonincreasing}
	The function $g_l$ is said to have a \emph{non-increasing lower tail} if there exists $x'\in \mathcal{I}$ such that: 
	\begin{enumerate}
		\item $g_l(x)\leq g_l(x')$, for all elements of the set $\left\{ x\in \mathcal{I}|\ x>x'\right\}$, provided it is not empty,
		\item $g_l(x_1)\geq g_l(x_2)$, for all $x_1, x_2 \in \mathcal{I}$ with $x_1 \leq x_2 \leq x'$.
	\end{enumerate}
	The function $g_l$ is said to have a \emph{non-increasing lower tail with threshold} $d_l=\sup \mathcal{G}_l$, , where $\mathcal{G}_l$ contains all $x'$ satisfying Conditions (1) and (2).
\end{definition} 

Before we proceed, let us make two remarks regarding monotone tail functions.

\begin{remark}
	A non-decreasing (resp. non-increasing) function has both a non-decreasing (resp. non-increasing) lower tail, and a non-decreasing (resp. non-increasing) upper tail.
\end{remark}

\begin{remark}\label{Rem33}
	If a function has a monotone lower tail (resp. upper tail), then it can only have a monotone upper tail (resp. lower tail).
\end{remark}
In order to illustrate Remark \ref{Rem33}, consider a function $h$ having a non-decreasing lower tail. This means that $h$ is a non-decreasing injection on Quadrant I of Figure \ref{Fig:Tails1a}, and cannot take values in Quadrant II or Quadrant III. Thus, $h$ can only take values in Quadrant IV, implying that it can only have a monotone upper tail.

The following lemma provides a link between the four different forms of monotone tail functions, and allows to focus on one of them in the proof of the main result. Throughout, for any set $\mathcal{D}$, the notation $\bar{\mathcal{D}}$ refers to the transformation of the elements of $\mathcal{D}$ such that:
\begin{equation}\label{Negset}x\in\mathcal{D} \Longleftrightarrow -x \in\bar{\mathcal{D}}.\end{equation}

\begin{lemma}\label{Lemma-Link}
	Consider the function $h_u:\mathcal{I}\mapsto \mathcal{J}$ which has a non-decreasing upper tail with threshold $c_u$, and the sets $\bar{\mathcal{J}}$ and $\bar{\mathcal{I}}$ defined according to \eqref{Negset}. The following holds:
	\begin{itemize}
		\item The function $g_u:\mathcal{I}\mapsto \bar{\mathcal{J}}$, such that $g_u(x)=-h_u(x)$ has a non-increasing upper tail with threshold $c_u$.
		\item The function $h_l:\bar{\mathcal{I}}\mapsto \bar{\mathcal{J}}$, such that $h_l(x)=-h_u(-x)$ has a non-decreasing lower tail with threshold $-c_u$.
		\item The function $g_l:\bar{\mathcal{I}}\mapsto \mathcal{J}$, such that $g_l(x) = h_u(-x)$ has a non-increasing lower tail with threshold $-c_u$.
	\end{itemize}
\end{lemma}

\subsection{Intermediary results}
In what follows, functions having non-decreasing upper tail are the main focus. Later, the links derived in Lemma \ref{Lemma-Link} are used to extend the results to all four types of monotone tail functions. The subscript ``$u$'' is momentarily dropped, such that a function having a non-decreasing upper tail with threshold $c$ is denoted by $h$.

In order to prove our main result, the function $h$ is transformed into a non-decreasing function $\tilde{h}$. This transformed function is such that it coincides with the original function $h$ in the upper part exceeding $c$, whereas the lower part is constructed such that the function $\tilde{h}$ is non-decreasing. However, moving back from the setting of $\tilde{h}$ to that of $h$ raises some challenges, especially in case of discontinuity at $c$. Circumventing those challenges requires some preliminary steps, to which the present subsection is devoted.\\

Define the quantity $h^{\star}$ as follows:
\begin{equation}\label{Defhstar}
	h^{\star}=\underset{x<c}{\sup\ } h(x),
\end{equation}
The function $h$ only takes values exceeding $h^{\star}$ if $x\geq c$. The value associated with $c$ used to determine the quadrants is given by $h^{\star}$.

Define also the function $\tilde{h}:\mathcal{I}\mapsto \mathcal{J}$ as follows:
\begin{equation}\label{Def_function-htilde}
	\tilde{h}(x)
	=\left\{
	\begin{array}{ll}
		h^{\star}, & \text{ \qquad if } x< c,\\
		\max\left\{ h(c), h^{\star}\right\}, & \text{ \qquad if } x=c, \\
		h(x), & \text{ \qquad if } x>c.
	\end{array}%
	\right.
\end{equation}
The new function $\tilde{h}$ is non-decreasing on the whole interval $\mathcal{I}$.

Finally, define the probability $\pi_{c}$ as follows
$$\pi_{c}=\mathbb{P}\left[h\left(X\right)\leq h^{\star} \right].$$ 
Since $c$ is the infimum of the set $\mathcal{H}$, it follows from Condition (1) of Definition \ref{Def:Tailmonotonic} that $X< c$ implies $h(X)\leq h^{\star}$. However, $h(X)\leq h^{\star}$ does not necessarily imply that $X\leq c$, especially if $h^{\star}$ corresponds to a horizontal part. Thus, the following inequality holds:
$$ \pi_{c}\geq \mathbb{P}\left[X< c \right],$$
and becomes an equality if $h$ is strictly increasing and continuous in $c$.\\

Four important lemmas are now stated and will be used for the proof of the main theorem; see Appendices \ref{Proofh-eq11}, \ref{L-h-1-1-Proof}, \ref{L-h-1-2-Proof} and \ref{L-h-1-3-Proof} for the proofs. The first lemma states that transforming a random variable using $h$ or $\tilde{h}$ does not affect the resulting cdf after the threshold $h^\star$.

\begin{lemma}\label{L-h-1-0} Consider a function $h$ having a non-decreasing upper tail with threshold $c$, and the function $\tilde{h}$ defined in \eqref{Def_function-htilde}. For $\pi_{c}=\mathbb{P}\left[h\left(X\right)\leq h^{\star} \right]$ and $h^{\star} = \underset{x<c}{\sup}h(x)$, it follows that: 
	\begin{equation}\label{h-eq1} 
		\mathbb{P}\left[\tilde{h}(X)\leq y \right]=\mathbb{P}\left[h(X)\leq y \right], \qquad \text{for }  y\geq h^{\star},
	\end{equation}
	and in particular:
	\begin{equation}
		\mathbb{P}\left[\tilde{h}(X)\leq h^{\star} \right]= \pi_{c}.\label{h-eq3}
	\end{equation} 
\end{lemma}

The second lemma states that the generalized quantile of the random variables $h(X)$ and $\tilde{h}(X)$ coincide when the probability level $p$ is taken sufficiently large.

\begin{lemma}\label{L-h-1-1} Consider a function $h$ having a non-decreasing upper tail with threshold $c$, and the function $\tilde{h}$ defined in \eqref{Def_function-htilde}. For $\pi_{c}=\mathbb{P}\left[h\left(X\right)\leq h^{\star} \right]$ and $h^{\star}=\sup_{x<c} h(x)$, it follows that: 
	\begin{eqnarray}
		F_{\tilde{h}(X)}^{-1(\alpha)}(p)&=&F_{h(X)}^{-1(\alpha)}(p),\qquad \text{for } p>\pi_{c} \text{ \ and \ } \alpha\in [0,1].  \label{Lh1-eq1}
	\end{eqnarray} 
\end{lemma}

The third lemma states that the same value is obtained regardless of whether the generalized quantile $F_X^{-1(\alpha)}(p)$ is transformed using $\tilde{h}$ or $h$, provided the value $p$ is sufficiently large.
\begin{lemma}\label{L-h-1-2}  Consider a function $h$ having a non-decreasing upper tail with threshold $c$, and the function $\tilde{h}$ defined in \eqref{Def_function-htilde}. For $\pi_{c}=\mathbb{P}\left[h\left(X\right)\leq h^{\star} \right]$, it follows that:
	\begin{equation}\label{Lh2-eq1}
		\tilde{h}\left(F_X^{-1(\alpha)}(p) \right)=h\left(F_X^{-1(\alpha)}(p) \right),\qquad \text{for } p>\pi_{c} \text{ \ and \ } \alpha\in [0,1].
	\end{equation}
\end{lemma}

The fourth lemma shows that left/right continuity of the function $h$ is inherited by $\tilde{h}.$

\begin{lemma}\label{L-h-1-3}  Consider a function $h$ having a non-decreasing upper tail with threshold $c$, and the function $\tilde{h}$ defined in \eqref{Def_function-htilde}. It follows that:
	\begin{enumerate}
		\item $h$ is left-continuous $\Rightarrow$ $\tilde{h}$ is left-continuous. 
		\item $h$ is right-continuous $\Rightarrow$ $\tilde{h}$ is right-continuous. 
	\end{enumerate}
\end{lemma}

\subsection{Main theorem}
Using these lemmas, together with Lemma \ref{Lemma-Link} that links the four types of monotone tail functions, it is now possible to derive conditions under which the inverse of a transformed random variable $h(X)$ can be expressed as the function $h$ applied on the inverse of the initial random variable $X$. The results are summarized in Theorem \ref{Theorem-h-1} and the proof can be found in Appendix \ref{Theorem-h-1-Proof}.

\begin{theorem} \label{Theorem-h-1}
	Consider a random variable $X$ with cdf $F_X$.
	
	\begin{itemize}
		\item For a function $h_u$ having a non-decreasing upper tail with threshold $c_u$, with $\pi_{c_u}=\mathbb{P}\left[h_u\left(X\right)\leq h^{\star}_u \right]$, and $h^{\star}_u = \underset{x<c_u}{\sup}h_u(x)$, the following holds:
		\begin{enumerate}
			\item If $h_u$ is left-continuous, then: 
			\begin{equation}\label{T-h-11}
				F_{h_u(X)}^{-1}(p)=h_u\left(F_X^{-1}(p) \right), \text{\qquad for} \ p>\pi_{c_u}.
			\end{equation}
			\item If $h_u$ is right-continuous, then:
			\begin{equation}\label{T-h-12}
				F_{h_u(X)}^{-1+}(p)=h_u\left(F_X^{-1+}(p) \right), \text{\qquad for} \ p
				\geq \pi_{c_u}.
			\end{equation}
		\end{enumerate}
		
		\item For a function $g_u$ having a non-increasing upper tail with threshold $d_u$, with $\lambda_{d_u}=\mathbb{P}\left[g_u\left(X\right)\geq g^{\star}_u \right]$, and $g^{\star}_u = \underset{x<d_u}{\inf}g_u(x)$, the following holds:
		\begin{enumerate}
			\item If $g_u$ is left-continuous, then: 
			\begin{equation}\label{T-h-21}
				F_{g_u(X)}^{-1+}(p)=g_u\left(F_X^{-1}(1-p) \right), \text{\qquad for} \ p<1-\lambda_{d_u}.
			\end{equation}
			\item If $g_u$ is right-continuous, then:
			\begin{equation}\label{T-h-22}
				F_{g_u(X)}^{-1}(p)=g_u\left(F_X^{-1+}(1-p) \right), \text{\qquad for} \ p
				\leq 1-\lambda_{d_u}.
			\end{equation}
		\end{enumerate}
		
		\item For a function $h_l$ having a non-decreasing lower tail with threshold $c_l$, with $\lambda_{c_l}=\mathbb{P}\left[h_l\left(X\right)\geq h^{\star}_l \right]$, and $h^{\star}_l = \underset{x>c_l}{\inf}h_l(x)$, the following holds:
		\begin{enumerate}
			\item If $h_l$ is left-continuous, then: 
			\begin{equation}\label{T-h-31}
				F_{h_l(X)}^{-1}(p)=h_l\left(F_X^{-1}(p) \right), \text{\qquad for} \ p\leq 1 - \lambda_{c_l}.
			\end{equation}
			\item If $h_l$ is right-continuous, then:
			\begin{equation}\label{T-h-32}
				F_{h_l(X)}^{-1+}(p)=h_l\left(F_X^{-1+}(p) \right), \text{\qquad for} \ p
				< 1- \lambda_{c_l}.
			\end{equation}
		\end{enumerate}
		
		\item For a function $g_l$ having a non-increasing lower tail with threshold $d_l$, with $\pi_{d_l}=\mathbb{P}\left[g_l\left(X\right)\leq g^{\star}_l \right]$, and $g^{\star}_l = \underset{x>d_l}{\sup}g_l(x)$, the following holds:
		\begin{enumerate}
			\item If $g_l$ is left-continuous, then: 
			\begin{equation}\label{T-h-41}
				F_{g_l(X)}^{-1+}(p)=g_l\left(F_X^{-1}(1- p) \right), \text{\qquad for} \ p\geq \pi_{d_l}.
			\end{equation}
			\item If $g_l$ is right-continuous, then:
			\begin{equation}\label{T-h-42}
				F_{g_l(X)}^{-1}(p)=g_l\left(F_X^{-1+}(1-p) \right), \text{\qquad for} \ p
				> \pi_{d_l}.
			\end{equation}
		\end{enumerate}
	\end{itemize}

\end{theorem}

Note that Expression \eqref{T-h-11} holds for values $p$ which are greater or equal to $\pi_c$. Expression \eqref{T-h-12}, on the other hand, only holds for values $p$ which are strictly larger than $\pi_c$. It is possible to construct counter-examples to prove that in general:
$$F_{h_u(X)}^{-1}\left(\pi_{c_u} \right)  \neq h_u\left(F_X^{-1}\left(\pi_{c_u} \right) \right),$$
for a function $h_u$ having a non-decreasing upper tail with threshold $c_u$. This is the case for a random variable $X$ with probability mass function:
$$\mathbb{P}\left[X=k \right]=\frac{1}{4},\ \text{for}\ k=0,1,2,3$$
and a function $h_u: [0,3]\rightarrow \mathbb{R}$ given by:
$$
h_u(x)= \left\{
\begin{array}{ll}
	x, & \text{ \qquad if } x \in [0,1],\\
	2-x, & \text{ \qquad if } x\in (1,2], \\
	2(x-2) & \text{ \qquad if } x\in (2,3]).
\end{array}%
\right.
$$
This remark is valid for the analogous expressions for the other types of monotone tail functions.

\section{Examples and associated results}\label{Section-Results}
\subsection{Option trading strategies}
Option trading strategies provide numerous examples where the random variable denoting the payoff of the strategy can be expressed using the stock price transformed by a monotone tail function. Such strategies are combinations of calls and puts, and sometimes the stock price itself, and are used in certain situations to construct neutral portfolios \citep{Kolb}. 

A common example is the \textit{straddle} option strategy. It consists in purchasing simultaneously a call option and a put option on the same underlying stock price, with the same strike and maturity. A straddle strategy is profitable especially when the underlying stock price rises. It is then relevant to determine the quantiles of the payoff in the upper tail.

For a given strike $K$, the payoff function of the straddle strategy is given by:
$$h(x)=\max\{K-x;0\} + \max\{x-K;0\},$$
with $h:\mathbb{R}^{+}\mapsto \mathbb{R}^{+}$. The function $h$ has a non-decreasing upper tail with threshold $2K$. Moreover, $h^\star = K$. Therefore, for a stock price at maturity denoted by $X$, the payoff at maturity is given by $h(X)$, and the corresponding quantile at the level $p$ has the following expression:
$$F_{h(X)}^{-1}(p) = \max\left\{K - F^{-1}_X(p);0\right\} + \max\left\{F^{-1}_X(p)-K;0\right\}, \qquad \text{for } p> \mathbb{P}\left[h(X)\leq K\right].$$

Many other option strategies are characterized by a payoff function having a monotone tail. The list includes the \textit{strangle} strategy, whose difference with the straddle is that the strike of the call option is higher than that of the put option, The payoff function in this case again has a non-decreasing upper tail. The payoff of the \textit{back spread with calls} (selling a call with strike $K_1$ and buying two calls with strike $K_2>K_1$) also has a non-decreasing upper tail. In contrast, the short versions of the straddle and the strangle involve non-increasing upper tails. The \textit{back spread with puts} (buying two puts with strike $K_1$ and selling a put with strike $K_2>K_1$) is an example of strategies whose payoff function has a non-increasing lower tail.

\subsection{Endowment assurances}
Consider an endowment assurance contract, with a benefit amount $S_1>0$ payable upon death if it occurs before the term $n$, or a maturity benefit $S_2>0$ if the policyholder survives the term. Suppose that the death and survival benefits are not equal, i.e. $S_1 \neq S_2$, which is not an uncommon feature of endowment assurances, and is the case under double endowment assurances for instance. Let $X$ be the random complete remaining lifetime of the policyholder. Assuming a constant discounting factor $v\in(0,1)$ throughout the term, the random present value of this benefit is $h(X)$, with $h:\mathbb{R}^{+}\mapsto \mathbb{R}^{+}$, where:
$$h(x)=S_1 v^x \mathbb{I}[x<n] + S_2 v^n \mathbb{I}[x\geq n],$$
and $\mathbb{I}$ is the indicator function. Note that the function $h$ is constant for $x\geq n$. Moreover, $h$ is right continuous at $n$, and has a jump at $n$ for $S_1 \neq S_2$.

If $S_1 \leq S_2v^n$, the function $h$ has a non-decreasing upper tail with threshold $n$, and hence:
$$F_{h(X)}^{-1+}(p) = S_1 v^{F_{X}^{-1+}(p)} \mathbb{I}\left[F_{X}^{-1+}(p)<n\right] + S_2 v^n \mathbb{I}\left[F_{X}^{-1+}(p)\geq n\right], \qquad \text{for } p \geq \mathbb{P}\left[h(X) \leq S_1\right]. $$
If $S_2v^n<S_1$, the function $g$ has a non-increasing lower tail with threshold $\frac{\log(S_2v^n)-\log(S_1)}{\log(v)}$, and hence:
$$F_{h(X)}^{-1}(p) = S_1v^{F_{X}^{-1+}(1-p)}\mathbb{I}\left[F_{X}^{-1+}(1-p)<n\right] + S_2 v^n \mathbb{I}\left[F_{X}^{-1+}(1-p)\geq n\right], \qquad \text{for } p > \mathbb{P}\left[h(X) \leq S_2v^n\right]. $$

\section{Quadrant perfect dependence}\label{SectionQuadrantDependence}
\cite{Cheung2009} introduced the concept of upper comonotonicity to model random variables which are perfectly positively dependent beyond a threshold of the joint support. More specifically, a random vector which is comonotonic in the top-right quadrant of the joint support is said to be upper comonotonic. Similar to comonotonic vectors, an upper-upper comonotonic vector has desirable additive properties for the Value-at-Risk of the sum of its components (see Proposition 6 in  \cite{Cheung2009}). In this section, we investigate how this concept integrates within the framework of monotone tail functions.

The four different types of monotone tail functions call for the introduction of three other dependence structures. These dependence structures were mentioned in e.g. \cite{DongCheungYang}, and we label them  \textit{lower-lower comonotonicity} (i.e.\  perfect positive dependence in the bottom-left quadrant), \textit{upper-lower counter-monotonicity} (i.e.\ perfect negative dependence in the bottom-right quadrant), and \textit{lower-upper counter-monotonicity} (i.e.\ perfect negative dependence in the top-left quadrant). In order to avoid confusion, we use the term \textit{upper-upper comonotonicity} when we refer to \cite{Cheung2009}'s dependence structure. The choice of this terminology is motivated by the following reasoning. For instance, an upper-lower counter-monotonic vector $\left(X_1,X_2\right)$ is counter-monotonic in the quadrant where $X_1$ is in the upper part of its univariate support and $X_2$ is in the lower part of its univariate support.

Throughout this section, we use the following notations for the four quadrants at $\underline{a}=\left(a_1,a_2\right)\in \mathbb{R}^2$, where the boundary is excluded:
\begin{eqnarray*}
	BL\left( \underline{a}\right) & = & (-\infty, a_1) \times (-\infty, a_2),\\
	BR\left( \underline{a}\right) & = & (a_1,+\infty) \times (-\infty, a_2),\\
	TL\left( \underline{a}\right) & = & (-\infty, a_1) \times (a_1,+\infty),\\
	TR\left( \underline{a}\right) & = & (a_1,+\infty) \times (a_2,+\infty).
\end{eqnarray*}
In particular, $BL(\underline{a})$ is the bottom-left quadrant, $BR(\underline{a})$ is the bottom-right quadrant, $TL(\underline{a})$ is the top-left quadrant, and $TR(\underline{a})$ is the top-right quadrant. We also use the bar to indicate that the boundary is included:
\begin{eqnarray*}
	\overline{BL}\left( \underline{a}\right) & = & (-\infty, a_1] \times (-\infty, a_2],\\
	\overline{BR}\left( \underline{a}\right) & = & [a_1,+\infty) \times (-\infty, a_2],\\
	\overline{TL}\left( \underline{a}\right) & = & (-\infty, a_1] \times [a_2,+\infty),\\
	\overline{TR}\left( \underline{a}\right) & = & [a_1,+\infty) \times [a_2,+\infty).
\end{eqnarray*}

\begin{definition}[Upper-upper comonotonicity]\label{Def:UpperCom}
	A two-dimensional random vector $\underline{X}$ with support $\mathcal{S}$ is upper-upper comonotonic with threshold $\underline{a}$ if the following conditions hold:
	\begin{enumerate}
		\item $\mathcal{S}_{\underline{a}}=\mathcal{S}\cap TR(\underline{a})$ is a comonotonic set, 
		\item $\mathbb{P}\left[\underline{X}>\underline{a} \right]>0$,
		\item $\mathcal{S}\cap \left(\mathbb{R}^2\setminus \left(TR(\underline{a})\cup \overline{BL}(\underline{a})\right)\right)$ is empty.
	\end{enumerate}
\end{definition}


\begin{definition}[Lower-Lower comonotonicitiy]\label{Def:LowerCom}
	A two-dimensional random vector $\underline{X}$ with support $\mathcal{S}$ is lower-lower comonotonic with threshold $\underline{a}$ if the following conditions hold:
	\begin{enumerate}
		\item $\mathcal{S}_{\underline{a}}=\mathcal{S}\cap BL(\underline{a})$ is a comonotonic set,
		\item $\mathbb{P}\left[\underline{X} <\underline{a}\right]>0$,
		\item $\mathcal{S}\cap \left(\mathbb{R}^2\setminus \left(BR(\underline{a})\cup \overline{TR}(\underline{a})\right)\right)$ is empty.
	\end{enumerate}
\end{definition}

\begin{definition}[Upper-Lower counter-monotonicity]\label{Def:UpperLowerCom}
	A two-dimensional random vector $\underline{X}$ with support $\mathcal{S}$ is upper-lower counter-monotonic with threshold $\underline{a}$ if the following conditions hold:
	\begin{enumerate}
		\item $\mathcal{S}_{\underline{a}}=\mathcal{S}\cap BR(\underline{a})$ is a counter-monotonic set,
		\item $\mathbb{P}\left[X_1>a_1,X_2<a_2\right]>0$,
		\item $\mathcal{S}\cap \left(\mathbb{R}^2\setminus\left(BR(\underline{a})\cup\overline{TL}(\underline{a})\right)\right)$ is empty.
	\end{enumerate}
\end{definition}

\begin{definition}[Lower-Upper counter-monotonicity]\label{Def:LowerUpperCom}
	A two-dimensional random vector $\underline{X}$ with support $\mathcal{S}$ is lower-upper counter-monotonic with threshold $\underline{a}$ if the following conditions hold:
	\begin{enumerate}
		\item $\mathcal{S}_{\underline{a}}=\mathcal{S}\cap TL(\underline{a})$ is a counter-monotonic set,
		\item $\mathbb{P}\left[X_1<a_1,X_2>a_2\right]>0$,
		\item $\mathcal{S}\cap \left(\mathbb{R}^2\setminus \left(TL(\underline{a})\cup \overline{BR}(\underline{a})\right)\right)$ is empty.
	\end{enumerate}
\end{definition}

The following lemma provides a link between the four different types of quadrant perfect dependence. A proof of the first statement can be found in Appendix \ref{Proof:LemmaLink2}, and similar arguments based on the rotation of the support of $(X_1,X_2)$ can be used to prove the remaining statements.

\begin{lemma}\label{LemmaLink2}
	Consider the random vector $\underline{X}=\left(X_1,X_2\right)$ with support $\mathcal{S}$ which is upper-upper comonotonic with threshold $(a_1,a_2)\in \mathcal{S}$. The following holds:
	\begin{itemize}
		\item The random vector $(-X_1,-X_2)$ is lower-lower comonotonic with threshold $(-a_1,-a_2)$.
		\item The random vector $(X_1,-X_2)$ is upper-lower counter-monotonic with threshold $(a_1,-a_2)$.
		\item The random vector $(-X_1,X_2)$ is lower-upper counter-monotonic with threshold $(-a_1,a_2)$.
	\end{itemize}
\end{lemma}

The following theorem establishes a relationship between quadrant perfect dependence and functions having a monotone tail. A proof of the first statement can be found in Appendix \ref{Proof:L-1}, together with a proof of how the second statement follows from the first (in combination with Lemmas \ref{Lemma-Link} and \ref{LemmaLink2}). A similar approach can be used to prove the remaining statements.
\begin{theorem}\label{L-1}
	Consider the comonotonic random vector $(X_1,X_2)$. Define the function $h$ as follows: 
	\begin{equation}\label{H_Alpha}h_\alpha(p) =F_{X_1}^{-1(\alpha)}(p) -F_{X_2}^{-1(\alpha)}(p),\qquad \text{for}\ p\in (0,1),\end{equation}
	with $\alpha \in [0,1]$, and for $\pi \in (0,1)$, we define the point $\underline{a}$ as follows:
	\begin{equation}\label{UpperComonotonic}
	\underline{a}=\left( F_{X_2}^{-1(\alpha)}(\pi),h_\alpha(\pi)\right).
	\end{equation}
	The following holds:
	\begin{itemize}
		\item $(X_2, X_1-X_2)$ is upper-upper comonotonic with threshold $\underline{a}$ iff the function $h_\alpha$ has a non-decreasing upper tail with threshold $\pi$.
		\item $(X_2, X_1-X_2)$ is lower-lower comonotonic with threshold $\underline{a}$ iff the function $h_\alpha$ has a non-decreasing lower tail with threshold $\pi$.
		\item $(X_2, X_1-X_2)$ is upper-lower counter-monotonic with threshold $\underline{a}$ iff the function $h_\alpha$ has a non-increasing upper tail with threshold $\pi$.
		\item $(X_2, X_1-X_2)$ is lower-upper counter-monotonic with threshold $\underline{a}$ iff the function $h_\alpha$ has a non-increasing lower tail with threshold $\pi$.
	\end{itemize}
\end{theorem}

Note that the link with quadrant perfect dependence in the above theorem is derived by assuming that $X_1$ and $X_2$ are comonotonic. Investigating whether this link is valid for other dependence structures of $\left(X_1,X_2\right)$ may deserve a separate contribution.

\section{Capital reducers}\label{Section-CapitalReducers}

\subsection{General setting}\label{Subsection21}
Consider an insurance company facing a random loss whose present value is denoted by $R_1$. The insurance company is required to hold a certain amount of cash such that the liabilities can be covered in the future with a sufficiently high probability. This capital buffer, or solvency capital, is determined using a risk measure, denoted by $\varphi$ \citep{Haezendonck_Goovaerts:1982,Wirch_and_Hardy2000,Hamidi2014,Boonen17,Happersberger2020,Chen2022}. 

In order to reduce the risk associated with its liabilities, and hence to decrease the amount of required regulatory capital, the insurance company can add an adequate partial hedge to the initial position $R_1$. Assume that today is time $0$, and let $R_2$ denote the present value of the future time-$T$ cash flow received by the insurer of a candidate hedger. The asset $R_2$ is expected to partially offset the losses from the liability $R_1$, and we assume that it can be acquired at time $0$ at a predetermined price $\varrho[R_2]$. Note that $R_2$ is an income random variable.

The present value of the insurer's liability after the income received from the asset $R_2$ is taken into account, is denoted by the loss random variable $Z$:
$$Z=R_1-R_2.$$
Positive values of $Z$ have to be interpreted as losses for the insurance company. If $R_2$ is a perfect hedge for $R_1$, we have that $Z=0$, almost surely. However, in a realistic situation, $R_1$ has a complex structure for which a perfect hedge does not exist. This means that $Z$ can still take positive values (i.e.\ $\mathbb{P}[Z>0]$ is positive), indicating a loss for the insurance company. Hence, $Z$ is the residual liability of the insurance company, i.e.\ the part of the losses of $R_1$ that cannot be countered by gains from $R_2$. The new position of the company now consists of the residual liability $Z$ and the price $\varrho[R_2]$ paid for acquiring $R_2$:
$$Z+\varrho[R_2].$$

Intuitively, $R_2$ is expected to be an adequate partial hedge if it is strongly positively dependent with $R_1$. We consider in this section the case where the random variables $R_1$ and $R_2$ are comonotonic, and we can write:
\begin{equation}\label{Eq1-1}
(R_1,R_2)\overset{d}{=} \left( F_{R_1}^{-1}(U), F_{R_2}^{-1}(U)\right),
\end{equation}
where $U$ is uniformly distributed over $[0,1]$. A formal definition of comonotonicity as well as some associated results and applications can be found in \cite{Dhaene2002a,Dhaene2002b}.

\subsection{Capital reducers}
An asset $R_2$ is said to be an adequate risk reducer for the liability $R_1$ if the residual liability $Z=R_1-R_2$ is `more preferable' than the initial liability $R_1$. The notion of `more preferable' can be quantified by using either integral stochastic orders \citep{Denuit:Modelling_Dependent_Risks,ShakedShantikumar}, or relevant risk measures. For studies on ordering of sums of random vectors, we refer to \cite{DenuitGenestMarceau} and \cite{CheungDhaeneKukushLinders}. We focus here on the latter, and we quantify the situation where $R_2$ is an effective risk reducer for $R_1$ by using a translation invariant risk measure $\varphi$. In particular, $R_2$ is a risk reducer for $R_1$ if the following inequality holds:
\begin{equation}\label{CapitalReducer-1}
\varphi\left[Z\right]+ \varrho[R_2]\leq \varphi[R_1].
\end{equation}
If \eqref{CapitalReducer-1} holds, we say that $R_2$ is a $\varphi$-reducer for the liability $R_1$.

Inequality \eqref{CapitalReducer-1} has a meaningful interpretation in terms of solvency capital reduction. Suppose that the amount of solvency capital a company has to hold is determined by using some risk measure $\varphi$. Then, the right-hand side of Inequality \eqref{CapitalReducer-1} corresponds with the amount of capital the insurance company has to hold to cover losses coming from the initial liability $R_1$, i.e.\ the liability the insurance company faces before the risk reducer is applied. The left-hand side of Inequality  \eqref{CapitalReducer-1} corresponds with the total amount of capital needed to buy the risk reducer and to set up the solvency capital one has to hold to cover the remaining losses after the income from the risk reducer has been taken into account. Therefore, Inequality \eqref{CapitalReducer-1} implies that buying the risk reducer $R_2$ leads to a reduction in necessary capital and is therefore considered to be the preferred strategy.

 This idea of risk reducers was introduced in \cite{CheungDhaeneLoTang} in the convex order sense, where the authors evaluate the effectiveness of a cash flow which is comonotonic with the initial liability\footnote{Note that the original definition is equivalently given for counter-monotonic risk reducers, where they consider the sum of random variables instead of their difference.}. A further generalization of convex risk reducers is studied in \cite{HeTangZhang_2016}. This is equivalent with setting $\varrho\left[R_2\right]=\mathbb{E}\left[R_2\right]$ in \eqref{CapitalReducer-1}, and setting $\varphi$ equal to TVaR$_p$ for all $p\in[0,1]$.
 
A definition in terms of TVaR$_p$ for all levels $p\in[0,1]$ is meaningful since it relates to the variability of the residual risk. Nevertheless, in practice, selecting the risk reducer $R_2$ using the TVaR for all levels $p$ can be too conservative, and insurance companies subject to VaR-based regulation might rather focus on the amount of capital reduction only. This is particularly relevant since the Solvency II directive prescribes the use of the VaR. Thus, in this section, we work in a simplified setting where we focus on one single level $p$, instead of all $p\in[0,1]$. In the sequel, we use the following definitions:
 
\begin{definition}\label{Def:VaRReducer}
	The cash flow $R_2$ is said to be a \emph{VaR-reducer} at the level $p$ for $R_1$ if 
	\begin{equation}\label{ConditionVaR}\varrho[R_2]\leq \emph{VaR}_p[R_1]-\emph{VaR}_p[Z],
	\end{equation}
	where $p\in [0,1]$ and $\varrho$ is a premium principle.
\end{definition}

Intuitively, Condition \ref{ConditionVaR} implies that the price $\varrho[R_2]$ we pay to set up the hedge should not be larger than the capital reduction we can obtain.
It is worth noting that if $R_2$ is a VaR-reducer at a level $p$ as well as for all $q\geq p$, then it is also a TVaR-reducer at the level $p$. In particular, by definition, $R_2$ is a VaR-reducer at a level $q$ if:
$$\text{VaR}_q[Z] + \varrho[R_2] \leq \text{VaR}_q[R_1].$$
If this inequality holds for all $q\geq p$, then:
$$\frac{1}{1-p}\int_p^1\left(\text{VaR}_q[Z] + \varrho[R_2]\right)\text{d}q\leq \frac{1}{1-p}\int_p^1\text{VaR}_q[R_1]\text{d}q,$$
and hence:
$$\text{TVaR}_q[Z] + \varrho[R_2] \leq \text{TVaR}_q[R_1].$$
That is, $R_2$ is a TVaR-reducer at a level $p$.

For $\varrho[R_2]=\mathbb{E}[R_2]$, \cite{CheungDhaeneLoTang} show that $R_2$ is a comonotonic TVaR-reducer if the vector $\left(R_2,R_1-R_2\right)$ is comonotonic (Theorem 1 and Lemma 6). Moreover, they provide examples where the intuition of merging $R_1$ and $-R_2$ to reduce the variance does not hold. Thus, Condition \ref{ConditionVaR} is not trivial. In what follows, we derive sufficient conditions for $R_2$ to be a VaR-reducer. Note that the setting considered in this section is similar to that leading to the equivalence between quadrant perfect dependence and functions having monotone tails. Thus, the results of this section can be derived in either frameworks.

\subsection{Sufficient conditions on the price of VaR-reducers}
The random difference $Z=R_1-R_2$ can be expressed as $Z \overset{d}{=} h(U)$, where:
	$$h(p)=F^{-1}_{R_1}(p)-F^{-1}_{R_2}(p),$$
	and $U$ is uniformly distributed over the unit interval. Taking into account Theorem \ref{Theorem-h-1}, we find that there are two conditions for $R_2$ to be a VaR-reducer.
	
The first condition follows from the case where $h$ has a non-decreasing upper/lower tail with threshold $c$. We find that $R_2$ is a VaR-reducer at the level $p$ for $R_1$ if its price $\varrho[R_2]$ satisfies:
	\begin{equation}\label{Reducer-Cond1}\varrho[R_2]\leq \text{VaR}_p[R_2],\end{equation}
	for $p>\pi_c$ when $h$ has a non-decreasing upper tail, or for $p\leq 1-\lambda_c$ when $h$ has a non-decreasing lower tail. Note that when $R_1$ is a loss random variable, we are mostly interested in large values for $p$, meaning that a function $h$ with non-decreasing upper tail is more useful.
	
Condition \eqref{Reducer-Cond1} provides a simple rule for determining a VaR-reducer, and it has a straightforward interpretation when combined with the concepts of upper-upper comonotonicity and lower-lower comonotonicity. In particular, we have shown in Theorem \ref{L-1} that upper-upper comonotonicity can be linked to the fact that the function $h$ has a non-decreasing upper tail, whereas lower-lower comonotonicity can be linked to the fact that $h$ has non-decreasing lower tail. Suppose that the losses in $R_1$ are large, i.e.\ $R_1$ takes large positive values. Then comonotonicity of $(R_1,R_2)$ implies that $R_2$ will be large too. Assume that $R_2$ is exceeding the comonotonic threshold. Then upper-upper comonotonicity of the vector $(R_2, R_1-R_2)$ ensures that $R_2$ will not dominate $R_1$, which would result in a heavy lower tail. Indeed, the loss  $R_1-R_2$ increases if $R_2$ increases. A similar rationale holds when $R_2$ is below a certain threshold, where lower-lower comonotonicity comes into play.

The second condition follows from the case where $h$ has a non-increasing upper/lower tail with threshold $c$ and $R_1$ and $R_2$ are continuous, we find that $R_2$ is a VaR-reducer at the level $p$ for $R_1$ if its price $\varrho[R_2]$ satisfies:
\begin{equation}\label{Reducer-Cond2}\varrho[R_2]\leq \text{VaR}_{p}[R_1]-\text{VaR}_{1-p}[R_1] + \text{VaR}_{1-p}[R_2],\end{equation}
for $p<1-\lambda_c$ when $h$ has a non-increasing upper tail, or for $p\geq \pi_c$ when $h$ has a non-increasing lower tail.

Condition \eqref{Reducer-Cond2} represents the opposite situation of Condition \eqref{Reducer-Cond1}. In particular, we know from Theorem \ref{L-1} that the properties of monotone tail functions leading to Inequality \eqref{Reducer-Cond2} correspond to quadrant counter-monotonicity. Thus, the loss $R_1 - R_2$ decreases if $R_2$ increases for the relevant values of $p$. This means that the cash flow $R_2$ dominates the initial position $R_1$ for those values.

Whereas Condition \eqref{Reducer-Cond1} can be satisfied as long as the hedger $R_2$ has a reasonable price, it is not guaranteed that the price of the hedger would satisfy Condition \eqref{Reducer-Cond2}. Everything held equal, this latter condition would be more likely satisfied in three situations: when the initial position $R_1$ has a heavy right tail (i.e. high extreme losses have a sufficiently high likelihood), when the initial position has a light left tail (i.e. very low losses have a sufficiently high likelihood), and when the candidate hedger has a heavy left tail (i.e. very low gains have a sufficiently high likelihood). In these cases, the upper bound is larger, and hence, the price of the hedger is likely to satisfy the condition.

It is worth noting that a direct consequence of the decomposition of quantiles in the tails is that the TVaR can also be decomposed in terms of either univariate TVaR's or LTVaR. In case the function $h$ has a non-decreasing upper/lower tail, the TVaR of $Z$ is equal to the difference of the TVaR's at the level $p$. In contrast, if the function $h$ has non-increasing upper/lower tail, the TVaR of $Z$ is equal to the difference of the LTVaR's at the level $1-p$ . Note that if the underlying random variables are not continuous, the TVaR's and LTVaR's have to be considered carefully as averages of either the left or right quantiles, depending on the nature of $h$. Taking into account these decompositions, Conditions \ref{Reducer-Cond1} and \ref{Reducer-Cond2} can be extended to TVaR-reducers in a straightforward way.

In what follows, we further analyze these two conditions under some fairly general distributional assumptions.

\subsection{Analysis of the conditions for being a VaR-reducer}
\subsubsection{Distributional assumptions and monotonic properties}
We consider the continuous random variable $W$ defined on the real line, with cdf $F_{W}$. The inverse cdf is denoted by $F^{-1}_{W}$, with $F^{-1}_W(0)=-\infty$ and $F^{-1}_W(1)=+\infty$. The probability density function (pdf) is denoted by $f_{W}$ and is assumed to be defined on $\mathbb{R}$. Moreover, $W$ is assumed to be standardized, i.e. $\mathbb{E}[W]=0$ and $\text{Var}[W]=1.$ We assume in this section that the components of $(R_1,R_2)$ can be written as follows:
\begin{equation}\label{LN_marginals}
\log R_i=\mu_i + \sigma_iW,\qquad \text{for}\ i=1,2,
\end{equation}
where $\mu_i$ and $\sigma_i$ are finite, and $\sigma_i> 0$, $i=1,2$. The VaR of $R_i$ follows from Theorem \ref{Theorem-Dhaene}:
\begin{equation}\label{LN-inv}
\text{VaR}_p[R_i] = \exp\left(\mu_i + \sigma_i F^{-1}_{W}(p)\right),
\end{equation}
for $p\in [0,1]$ and $i=1,2$. 

The choice \eqref{LN_marginals} includes a large set of distributions. For example, one can assume $W$ to be a standard normal distribution, which implies that $R_1$ and $R_2$ follow a lognormal distribution. In order to allow more flexibility, one can also assume $W$ to come from a L\'evy distribution. The random variables $R_1$ and $R_2$ are then exponential L\'evy models. Such models have proven useful to model asset returns; see e.g.\ \cite{LevyProcesses_Book}. In \cite{VANBILSEN201919}, the authors use the model \eqref{LN_marginals}, where $W$ is assumed to be a standardized Variance Gamma distribution, to model stylized facts of asset returns such as heavy tails, leptokurticy and skewness. They show that this model is accurate for modeling stock returns, especially long term asset returns; see also \cite{Albrecher_et_al_2007} and \cite{Linders_et_al_2014}.

The function $h:\mathcal{I}\mapsto \mathcal{J}$ can be written as follows:
\begin{equation}\label{LNg}
h(u)=\text{e}^{\mu_1+\sigma_1F_{W}^{-1}(u)}-\text{e}^{\mu_2+\sigma_2F_{W}^{-1}(u)},
\end{equation}
where $\mathcal{I}=(0,1)$ and $\mathcal{J}=(y^{\min},y^{\max})$, with $y^{\min}=\underset{u\in\mathcal{I}}{\inf}\ h(u)$ and $y^{\max}=\underset{u\in\mathcal{I}}{\sup}\ h(u)$.

The following proposition provides conditions under which $h$ is monotone, and shows that when these conditions are not satisfied, then $h$ always has a monotone tail. We refer to Appendix \ref{LN-case2-proof} for a proof; see also \cite{Chaoubi} for a related result. Note that we discard the trivial case where $\sigma_1=\sigma_2$ and $\mu_1=\mu_2$, which leads to $Z=0$ almost surely.

\begin{proposition}\label{LN-case2-functiong} Consider the function $h$ defined in \eqref{LNg} and the comonotonic difference $Z\overset{d}{=}h(U)$.
	For $\sigma_1\neq \sigma_2$, define $c^{\star}$, $u^{\star}$ and $y^\star=h(u^{\star})$ as follows:
	\begin{equation}\label{Prop1-Parameters}
	c^{\star}=  F_W\left(\frac{\mu_2-\mu_1}{\sigma_1-\sigma_2} \right),\qquad \text{ and } \qquad
	u^{\star}= F_W\left(\frac{\log\left(\frac{\sigma_1}{\sigma_2}\right)+\mu_1-\mu_2}{\sigma_2-\sigma_1}\right).
	\end{equation}
	
	\begin{enumerate}
		\item If $\sigma_1=\sigma_2$ and $\mu_1>\mu_2$, we have that $y^{\min}=0$, $y^{\max}=+\infty$, and  $h$ is non-decreasing.
		\item If $\sigma_1=\sigma_2$ and $\mu_2>\mu_1$, we have that $y^{\min}=-\infty$, $y^{\max}=0$, and $h$ is non-increasing.
		\item If $\sigma_1>\sigma_2$, we have that  $y^{\min}=y^{\star}<0$, $y^{\max}=+\infty$, and the function $h$ has a non-decreasing upper tail with threshold $c^{\star}$.
		
		\item If $\sigma_1<\sigma_2$, we have that $y^{\min}=-\infty$,  $y^{\max}=y^{\star}>0$, and the function $h$ has a non-increasing upper tail with threshold $c^{\star}$.
	\end{enumerate}
\end{proposition}

We conclude from Proposition \ref{LN-case2-functiong} that there are four possible shapes for the function $h$, in case the marginal distributions are given by \eqref{LN_marginals}. When $\sigma_1=\sigma_2$, the function $h$ is always monotone. When $\sigma_1\neq \sigma_2$, the function $h$ is either convex or concave. Thus, $h$ can have a monotone upper tail, but cannot have a monotone lower tail.

\begin{figure}[!h]
	\centering
	\input{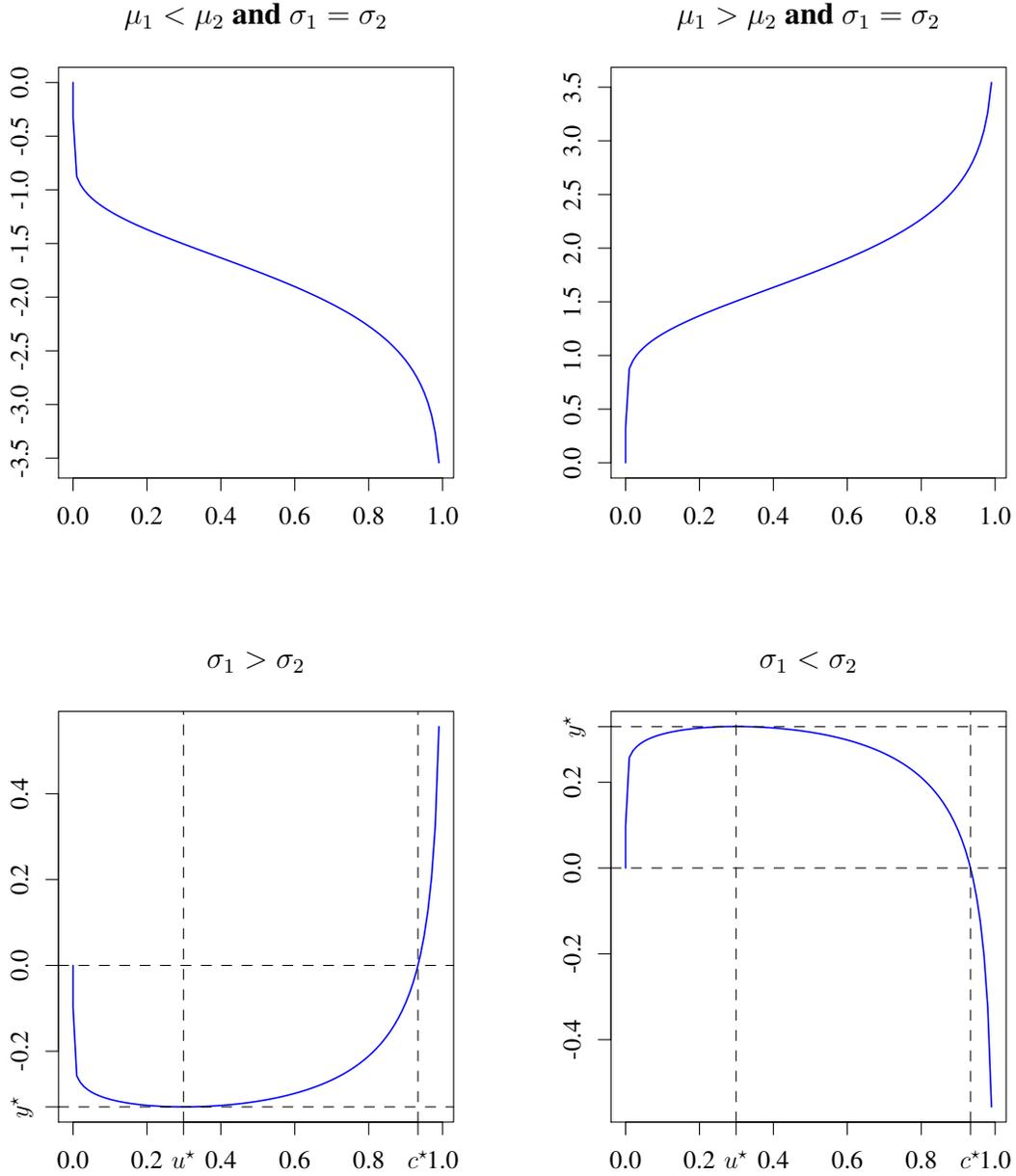}
	\caption{\footnotesize The function $h$ given by \eqref{LNg}. The top two panels are the situations where $\sigma_1=\sigma_2$, whereas the bottom two panels are the situations where $\sigma_1\neq \sigma_2$. The optimum $u^\star$ and the threshold $c^{\star}$ are highlighted with the vertical lines.}\label{Fig:LN-functiong}
\end{figure}

Figure \ref{Fig:LN-functiong} displays these four shapes. We assume in this figure that $W$ has a standard normal distribution, i.e.\ $R_1$ and $R_2$ are lognormal distributed. The top two panels correspond to the situations where $\sigma_1=\sigma_2$, i.e. the function $h$ is always monotone. The bottom two panels illustrate the case where $\sigma_1\neq \sigma_2$, with $h$ being either convex or concave. Thus, $h$ can have a monotone upper tail, but cannot have a monotone lower tail.

\subsubsection{Analysis}
Now, we further analyze Conditions \eqref{Reducer-Cond1} and \eqref{Reducer-Cond2} under the present distributional assumptions. First, notice that because $R_1$ and $R_2$ are continuous, the right and left inverses $F^{-1+}_Z$ and $F^{-1}_Z$ are equal. For $h^{\star}=h\left(c^{\star}\right)$, the uniform distribution of $U$ leads to
$$\pi_{c^{\star}}=c^{\star}, \qquad \text{ and } \qquad \lambda_{c^{\star}}= 1-c^{\star}.$$
 
Therefore, we conclude that:

\begin{enumerate}
	\item[Case 1:] If either $\sigma_1=\sigma_2$ and $\mu_1>\mu_2$, or $\sigma_1>\sigma_2$, then the hedger $R_2$ is a VaR-reducer at the level $p\geq c^{\star}$ for $R_1$ provided the price $\varrho[R_2]$ satisfies:
	\begin{equation}\label{eq-LN-1}
	\varrho[R_2]\leq \text{VaR}_p[R_2].
	\end{equation}
	\item[Case 2:] If either $\sigma_1=\sigma_2$ and $\mu_1<\mu_2$, or $\sigma_1<\sigma_2$, then the hedger $R_2$ is a VaR-reducer at the level $p\leq c^{\star}$ for $R_1$ provided the price $\varrho[R_2]$ satisfies:
	\begin{equation}\label{eq-LN-2}
	\varrho[R_2]\leq \text{VaR}_{1-p}[R_2] + \text{VaR}_p[R_1] - \text{VaR}_{1-p}[R_1]. 
	\end{equation}
\end{enumerate}

We first consider Case 1. If the standard deviation of the log-return of $R_2$ is smaller than that of the initial position, i.e.\ $\sigma_2<\sigma_1$, or if these standard deviations are equal but the expected log-return of $R_1$ is larger than that of $R_2$, then any hedger $R_2$ is a VaR-reducer provided the price $\varrho[R_2]$ is not `too high'. Indeed, it is reasonable to assume that any realistic price $\varrho[R_2]$ will be lower than the amount of capital required for holding $R_2$ only. Note that in this case, the comonotonic difference $Z$ is bounded from below by $y^{\star}$.

%
%


When $R_2$ satisfies the conditions of Case 1, its price must be such that $\varrho[R_2]\leq $VaR$_p[R_2]$ for all $p\geq c^{\star}$. Since VaR is non-decreasing with respect to $p$, then, taking into account \eqref{LN-inv} and \eqref{Prop1-Parameters}, we find that a sufficient condition on the price of $R_2$ under Case 1 is as follows:
$$\varrho\left[R_2\right]\leq\text{VaR}_{c^{\star}}[R_2]=\exp\left(\frac{\sigma_1\mu_2 - \sigma_2\mu_1}{\sigma_1-\sigma_2}\right).$$

Consider now Case 2, where $\sigma_2>\sigma_1$. In this case, reducing the risk of the liability $R_1$ by adding the cash flow $R_2$ `removes' the upper tail. Indeed, the random variable $R_2$ partially hedges $R_1$ and takes away part of the losses (i.e.\ the positive part) of $R_1$ and the random variable $Z$ is bounded from above. From Proposition \eqref{Prop1-Parameters}, one can verify that the maximum of the random variable $Z$ is given by:
\begin{equation}\label{NewPos-Upper} Z \leq  \text{e}^{\frac{\mu_1\sigma_2-\mu_2\sigma_1}{\sigma_2-\sigma_1}} 
\left(\text{e}^{\frac{\sigma_1 }{\sigma_2-\sigma_1}\log \frac{\sigma_1}{\sigma_2}}-\text{e}^{\frac{\sigma_2 }{\sigma_2-\sigma_1}\log \frac{\sigma_1}{\sigma_2}} \right),\end{equation}
whereas its minimum is $-\infty$. Thus, the losses are bounded but the gains are not. However, for $R_2$ to be a VaR-reducer, its price must satisfy Inequality \eqref{eq-LN-2}, which is not guaranteed. Therefore, there are some situations where the random variable $R_2$ is too expensive for the capital reduction it provides, and hence, reducing the liability $R_1$ by adding $R_2$ is not optimal..

Finally, note that Conditions \eqref{eq-LN-1} and \eqref{eq-LN-2} can be extended in a straightforward way to account for TVaR-reducers, and the above discussion remains relevant in this case. In particular, we have that:
\begin{enumerate}
	\item[Case 1:] If either $\sigma_1=\sigma_2$ and $\mu_1>\mu_2$, or $\sigma_1>\sigma_2$, then the hedger $R_2$ is a TVaR-reducer at the level $p\geq c^{\star}$ for $R_1$ provided the price $\varrho[R_2]$ satisfies:
	$$
	\varrho[R_2]\leq \text{TVaR}_p[R_2].
	$$
	\item[Case 2:] If either $\sigma_1=\sigma_2$ and $\mu_1<\mu_2$, or $\sigma_1<\sigma_2$, then the hedger $R_2$ is a TVaR-reducer at the level $p\leq c^{\star}$ for $R_1$ provided the price $\varrho[R_2]$ satisfies:
	$$
	\varrho[R_2]\leq \text{LTVaR}_{1-p}[R_2] + \text{TVaR}_p[R_1] - \text{LTVaR}_{1-p}[R_1].
	$$
\end{enumerate}

\section{Conclusion}\label{Section-Final}
We introduce the notion of monotone tail functions and analyze the properties of random variables transformed using such functions. We extend the results of \cite{Dhaene2002a} stating that the quantile of a random variable transformed with a monotone function is related to the transformed quantile of the original random variable. In particular, we show that this result holds for monotone tail functions, provided the quantile level is within an explicitly given interval. We also show that these functions are useful in different situations, such as the evaluation of the payoff of option trading strategies, or the present value of endowment assurances.

We address generalizations of the upper comonotonicity structure introduced in \cite{Cheung2009} through the lens of the proposed theory. We show how quadrant perfect dependence integrates within the notion of monotone tail functions in the framework of risk reducers.

We study the relevance of functions having monotone tails in the context of risk reduction, which is central to many financial and actuarial problems. Often, the problem boils down to selecting a hedger for an initial position, such that the level of the risk is reduced, and ultimately, the level of required regulatory capital too. We investigate conditions on the price of a comonotonic hedger ensuring a reduction of the Value-at-Risk, and we derive using the theory developed in this paper a simple rule for the choice of the comonotonic hedger.

\small
\bibliographystyle{agsm}
\bibliography{References_TailMonotonicity}
\appendix
\section{Appendix : proofs}\label{Section-Appendix}
\renewcommand{\theequation}{\thesection.\arabic{equation}}

\subsection{Proof of Lemma \ref{Lemma-Link}}\label{Proof-LemmaLink}
We prove the equivalence between $h_u$ having a non-decreasing upper tail and $g_u$ having a non-increasing upper tail. Similar arguments can be used to prove the remaining  statements.

If the function $h_u$ has a non-decreasing upper tail, there exists $x'\in\mathcal{I}$ such that Conditions (1) and (2) in Definition \ref{Def:Tailmonotonic} hold. On the one hand, Condition (1) states that for $x\in\mathcal{I}$ with $x<x'$, we have that $h_u(x)\leq h(x')$, which is equivalent to $-h_u(x)\geq -h_u(x')$. Thus, if we define the function $g_u:\mathcal{I}\mapsto \bar{\mathcal{J}}$ such that $g_u(x)=-h_u(x)$, we have that 
\begin{equation}\label{ProofLink-1}g_u(x) \geq g_u(x'), \qquad \text{for } x\in\mathcal{I} \text{ with } x<x'.\end{equation}
On the other hand, Condition (2) in Definition \ref{Def:Tailmonotonic} states that for any $x_1,x_2\in \mathcal{I}$ with $x'\leq x_1\leq x_2$, the inequality $h_u(x_1)\leq h_u(x_2)$ holds, which is equivalent with 
\begin{equation}\label{ProofLink-2}g_u(x_1)\geq g_u(x_2).\end{equation}
Combining \eqref{ProofLink-1} and \eqref{ProofLink-2} proves that $g_u$ satisfies the conditions of a non-increasing upper tail function from Definition \ref{Def:Tailmonotonic_nonincreasing}.

Further, we have from \eqref{ProofLink-1} and \eqref{ProofLink-2}, that any $x'$ satisfying the conditions from Definition \ref{Def:Tailmonotonic} of non-decreasing upper tail functions (i.e. $x'\in \mathcal{H}_u$), satisfies the conditions from Definition \ref{Def:Tailmonotonic_nonincreasing} of non-increasing upper tail functions. Thus, the infimum of set $\mathcal{H}_u$ associated with $h_u$ is equal to the infimum of set $\mathcal{G}_u$ associated with $-h_u$. Therefore, we conclude that $h_u$ and $g_u$ have the same threshold.

\subsection{Proof of Lemma \ref{L-h-1-0}}\label{Proofh-eq11}
We can use the law of total probability to write
	$$\mathbb{P}\left[\tilde{h}(X)\leq  y \right]=\mathbb{P}\left[X<c,\tilde{h}(X)\leq  y\right] + \mathbb{P}\left[X=c,\tilde{h}(X)\leq  y\right] + \mathbb{P}\left[X>c,\tilde{h}(X)\leq  y\right].$$
	If $X<c$, then $h(x) \leq h^{\star}$, and since $h^{\star}\leq y$, we find  
	$$\mathbb{P}\left[X<c,\tilde{h}(X)\leq  y\right]=\mathbb{P}\left[X<c\right].$$
For the second probability, if we use the definition  of $\tilde{h}$ and the fact that $h(c)\leq \max\{h(c),h^{\star}\}$ to write
	$$\mathbb{P}\left[X=c,\tilde{h}(X)\leq  y\right]=\mathbb{P}\left[X=c,h(c)\leq  y\right].$$		
	For the third probability, we have by definition that $\tilde{h}(x)=h(x)$ for $x>c$, and thus:
	$$\mathbb{P}\left[X>c,\tilde{h}(X)\leq  y\right]=\mathbb{P}\left[X>c,h(X)\leq  y \right].$$
	Combining the alternative representations of these three probabilities and using again the law of total probability, we find that $\mathbb{P}\left[\tilde{h}(X)\leq  y \right]=\mathbb{P}\left[h(X)\leq  y \right]$.

\subsection{Proof of Lemma \ref{L-h-1-1}} \label{L-h-1-1-Proof}
In order to prove \eqref{Lh1-eq1}, we start with proving the equality
\begin{equation}\label{LH-1}
F_{\tilde{h}(X)}^{-1}(p)=F_{h(X)}^{-1}(p),\qquad \text{for } p>\pi_{c},
\end{equation}
where $h$ has a non-decreasing upper tail. We have from \eqref{inverse} that:
$$F_{h(X)}^{-1}(p)=\inf\left\{y\in \mathbb{R}|\ \mathbb{P}\left[h(X)\leq y \right]\geq p \right\}	\overset{not.}{=}\inf \mathcal{A}_{h(X)},$$
$$	F_{\tilde{h}(X)}^{-1}(p)=\inf\left\{y\in \mathbb{R}|\ \mathbb{P}\left[\tilde{h}(X)\leq y \right]\geq p \right\}\overset{not.}{=}\inf \mathcal{A}_{\tilde{h}(X)}.$$
Since $p>\pi_{c}$, we have for $y\leq h^{\star}$ that $\mathbb{P}\left[ h(X)\leq y\right] <p$ from the definition of $\pi_c$, and that $\mathbb{P}\left[ \tilde{h}(X)\leq y\right] <p$ from \eqref{h-eq3}. This implies that $ y \notin \mathcal{A}_{h(X)}$ and $ y \notin \mathcal{A}_{\tilde{h}(X)}$ for $y\leq h^{\star}$, and hence, the sets $\mathcal{A}_{h(X)}$ and $\mathcal{A}_{\tilde{h}(X)}$ contain only values $y>h^{\star}$. It then follows from  \eqref{h-eq1} that $\mathcal{A}_{h(X)}=\mathcal{A}_{\tilde{h}(X)}$, which proves \eqref{LH-1}. 

We now prove 
\begin{equation}\label{LH-2}
F_{\tilde{h}(X)}^{-1+}(p)=F_{h(X)}^{-1+}(p),\qquad \text{for } p>\pi_{c}.
\end{equation}
We use \eqref{+inverse} to write: 
$$F_{h(X)}^{-1+}(p)=\sup\left\{y\in \mathbb{R}|\ \mathbb{P}\left[h(X)\leq y \right]\leq p \right\}	\overset{not.}{=} \sup \mathcal{B}_{h(X)},$$
$$	F_{\tilde{h}(X)}^{-1+}(p)=\sup\left\{y\in \mathbb{R}|\ \mathbb{P}\left[\tilde{h}(X)\leq y \right]\leq p \right\}\overset{not.}{=} \sup \mathcal{B}_{\tilde{h}(X)}.$$
Since $p>\pi_{c}$, we conclude that  $h^{\star}\in \mathcal{B}_{h(X)}$ using the definition of $\pi_c$. Similarly, we find that  $h^{\star} \in \mathcal{B}_{\tilde{h}(X)}$ when taking into account \eqref{h-eq3}. Thus, $\sup \mathcal{B}_{h(X)}\geq h^{\star}$  and $\sup \mathcal{B}_{\tilde{h}(X)}\geq h^{\star}$. The proof of \eqref{LH-2} is completed by using \eqref{h-eq1} to write: 
$$\sup \mathcal{B}_{h(X)}=	\sup\left\{y \geq h^{\star}|\ \mathbb{P}\left[\tilde{h}(X)\leq y \right]\leq p \right\}=\sup\mathcal{B}_{\tilde{h}(X)}.$$
Finally, combining \eqref{LH-1} and \eqref{LH-2} with the definition \eqref{alpha-inverse} proves \eqref{Lh1-eq1}. 

\subsection{Proof of Lemma \ref{L-h-1-2}}\label{L-h-1-2-Proof}
The inverse $F_X^{-1}$ is defined as follows:
$$F_X^{-1}(p)=\inf \left\{ y \in \mathbb{R}|\ F_X(y)\geq p \right\}=\inf \mathcal{A}_X.$$
It follows from the definition of $h^{\star}$ and $\pi_c$ that we can write:
\begin{equation}\label{L-h-1-2-Proof-1}
\pi_c=\mathbb{P}\left[X<c \right], \qquad \text{if} \ h(c)> h^{\star},
\end{equation}
and 
\begin{equation}\label{L-h-1-2-Proof-2}
\pi_c \geq \mathbb{P}\left[X\leq c \right], \qquad \text{if} \ h(c)\leq h^{\star}.
\end{equation}
We prove the result for these two cases separately.

\subsubsection*{Case 1: $h(c)>h^{\star}$}
Since $h(c)>h^{\star}$, we have that 
$$\tilde{h}(c)=\max \left\{h(c),h^{\star} \right\} = h(c).$$
which implies that:
\begin{equation}\label{L-h-1-2-Proof-3}
\tilde{h}(x)=h(x), \qquad \text{for} \ x\geq c.
\end{equation} 
Recall from \eqref{L-h-1-2-Proof-1} that in this case, we can write the probability $\pi_c$ as follows: $$\pi_c=\mathbb{P}\left[X<c \right],$$ 
which leads to:
\begin{equation}\label{L-h-1-2-Proof-5}
\mathbb{P}\left[X\leq c \right]=\pi_c + \mathbb{P}\left[X=c \right].
\end{equation}

If $\mathbb{P}\left[X=c \right]=0$, i.e.\ the cdf $F_X$ does not jump in $c$ and $\pi_c=F_X(c)$, we immediately find that $p>\pi_c$ implies
	$c\leq F_X^{-1}(p)$. If $\mathbb{P}\left[X=c \right]>0$, i.e.\ the cdf $F_X$ has a jump in $c$ with $F_X(c)>\pi_c$, then we have either $\pi_c <F_{X}(c)< p$ or $\pi_c<p \leq F_X(c)$. For $\pi_c <F_{X}(c)< p$, the element $p$ is beyond the jump part corresponding with $c$, and hence, the monotonicity of $F^{-1}_X$ implies that $c\leq F^{-1}_X(p)$. For $\pi_c<p \leq F_X(c)$, the element $p$ falls within the jump part, and hence, $c=F^{-1}_X(p)$. Therefore, we have in all cases that:
	\begin{equation}\label{L-h-1-2-Proof-5-2}
	c\leq F_X^{-1}(p), \qquad \text{for }p>\pi_c.
\end{equation}

We can use the definition of the $\alpha$-inverse to conclude that \eqref{L-h-1-2-Proof-5-2} implies 
$$c\leq F_X^{-1}(p) \leq F_X^{-1(\alpha)}(p).$$ 
Taking into account \eqref{L-h-1-2-Proof-3} proves that $h\left(F^{-1(\alpha)}_X(p)\right)=\tilde{h}\left(F^{-1(\alpha)}_X(p)\right)$ if $p>\pi_c.$ 

\subsubsection*{Case 2: $h(c)\leq h^{\star}$}
In case $h(c)\leq h^{\star}$, we find that the probability $\pi_c$ can be bounded from below as follows:
\begin{equation}\label{L-h-1-2-Proof-4}
\pi_c \geq F_X(c),
\end{equation}
which implies that for $p>\pi_c$, we have $p>F_X(c)$. Thus, from the definition of the inverse $F^{-1}_X$, we find that $c \notin \mathcal{A}_X$. Note that 
$$F_X\left(F^{-1}_X(p)\right)\geq p,$$
from which we find that $F_X^{-1}(p) \in \mathcal{A}_X$. 
Using the non-decreasingness of $F_X^{-1}$ together with the definition of the generalized $\alpha$-inverse, we have that:
$$c<F_X^{-1}(p)\leq F_X^{-1(\alpha)}(p).$$
Finally, from Expression \eqref{Def_function-htilde} of $\tilde{h}(x)$, we conclude that $h\left(F^{-1(\alpha)}_X(p)\right)=\tilde{h}\left(F^{-1(\alpha)}_X(p)\right)$. 

\subsection{Proof of Lemma \ref{L-h-1-3}}\label{L-h-1-3-Proof}
For $x<c$, $\tilde{h}(x)=h^{\star}$ and the function $\tilde{h}$ is continuous. For $x>c$, $\tilde{h}(x)=h(x)$, and the function $\tilde{h}$ has the same continuity as $h$. It is then sufficient to analyze the functions $\tilde{h}$ and $h$ in $c$.

Suppose that $h$ is left-continuous. The function $h$ has a non-decreasing upper tail with threshold $c$, from which we find that $h(c)\leq h^{\star}$, and hence $\tilde{h}(c)=h^{\star}$. Therefore, $\tilde{h}$ is left-continuous in $c$.

Suppose that $h$ is right-continuous. We have that $h(c)\geq h^{\star}$, which implies $\tilde{h}(c)=h(c)$, and therefore, $\tilde{h}(x)=h(x)$ for $x \geq c$. We can conclude that $\tilde{h}$ is right-continuous.

\subsection{Proof of Theorem \ref{Theorem-h-1}}\label{Theorem-h-1-Proof}
\subsubsection*{Expression \eqref{T-h-11}}
If the function $h_u$ is left-continuous with a non-decreasing upper tail, then $\tilde{h}_u$ as defined in \eqref{Def_function-htilde} is left-continuous and non-decreasing. Using \eqref{T-h-D1} gives:
\begin{equation}\label{Theorem-h1-eq1}
F_{\tilde{h}_u(X)}^{-1}(p)=\tilde{h}_u\left(F_X^{-1}(p) \right), \qquad \text{for } p\in (0,1).
\end{equation}
Using  Lemma \ref{L-h-1-1} with $\alpha=1$ gives
\begin{equation}\label{Theorem-h1-eq2}
F_{\tilde{h}_u(X)}^{-1}(p)=	F_{h_u(X)}^{-1}(p),\qquad \text{for } \pi_{c_u}<p,
\end{equation}
with $\pi_{c_u} = \mathbb{P}\left[h_u(X)\leq h^{\star}\right]$ and $h^\star = \underset{x<c_{u}}{\sup} h_u(x)$. If we put $\alpha=0$ in Lemma \ref{L-h-1-2}, we find
\begin{equation}\label{Theorem-h1-eq3}
\tilde{h}_u\left(F_X^{-1}(p) \right)=h_u\left(F_X^{-1}(p) \right),\qquad \text{for } \pi_{c_u}<p.
\end{equation}
Plugging Expressions \eqref{Theorem-h1-eq2} and \eqref{Theorem-h1-eq3} in \eqref{Theorem-h1-eq1} proves Expression \eqref{T-h-11}.

\subsubsection*{Expression \eqref{T-h-12}}
In order to prove Expression \eqref{T-h-12}, we start by noting that $\tilde{h}_u$ is now right-continuous and non-decreasing. Using Expression \eqref{T-h-D2} gives
\begin{equation}\label{Theorem-h1-eq4}
F_{\tilde{h}_u(X)}^{-1+}(p)=\tilde{h}_u\left(F_X^{-1+}(p) \right), \qquad \text{for } p\in (0,1).
\end{equation}
The proof for $p>\pi_{c_u}$ follows from Lemma \ref{L-h-1-1} and Lemma \ref{L-h-1-2}, with $\alpha=0$.

Let us now consider the case where $h_u$ is right-continuous with $p=\pi_{c_u}$. First, note that the right-inverse $F_{h_u(X)}^{-1+}$ is a right-continuous function, and hence:
$$F_{h_u(X)}^{-1+}(\pi_{c_u})=\lim_{p\downarrow \pi_{c_u}}F_{h_u(X)}^{-1+}(p).$$
Using the equality in \eqref{T-h-12} for the case $p>\pi_{c_u}$, we find:
$$F_{h_u(X)}^{-1+}(\pi_{c_u})=\lim_{p\downarrow \pi_{c_u}}h_u\left(F_{X}^{-1+}(p) \right).$$
Since $F_{X}^{-1+}$ is a non-decreasing function, the r.h.s.\ of this expression can be rewritten as
$$\lim_{p\downarrow \pi_{c_u}}h_u\left(F_{X}^{-1+}(p) \right)=\lim_{x\downarrow F_{X}^{-1+}(\pi_{c_u})}h_u\left(x \right).$$
Using the right continuity of the function $h_u$, leads to:
$$F_{h_u(X)}^{-1+}(\pi_{c_u})=h_u\left(F_{X}^{-1+}(\pi_{c_u}) \right),$$
which therefore proves the result for $p\geq \pi_{c_u}$.

\subsubsection*{Expressions \eqref{T-h-21} and \eqref{T-h-22}}
We know from Lemma \ref{Lemma-Link} that $h_u$ having a non-decreasing upper tail with threshold $c_u$ is equivalent with $g_u=-h_u$ having a non-increasing upper tail with threshold $c_u$. Moreover, the left or right continuity of $h_u$ is the same as that of $g_u$. We also have that $g^{\star}_u = -h^{\star}_u$, which leads to $\mathbb{P}\left[h_u(X)\leq h^{\star}_u\right]=\mathbb{P}\left[g_u(X)\geq g^{\star}_u\right]$.

Therefore, since $g_u = -h_u$, we find that Expression \eqref{T-h-11} is equivalent with:
$$F^{-1}_{-g_u(X)}(p) = -g_u\left(F^{-1}_X(p)\right), \qquad \text{for } p > \lambda_{c_u},$$
and using \eqref{T-h-D3} (or, equivalently, from \eqref{T-h-D4}), we have that $F_{-g_u(X)}^{-1}(p)=-F^{-1+}_{g_u(X)}(1-p)$, and rearranging for $1-p$ proves \eqref{T-h-21}.

Further, Expression \eqref{T-h-12} is equivalent with:
$$F_{-g_{u}(X)}^{-1+}(p) = -g_u\left(F^{-1+}_{X}(p)\right), \qquad \text{for } p\geq \lambda_{c_u},$$
and using again \eqref{T-h-D3} and rearranging for $1-p$  proves \eqref{T-h-22}.

Finally, using the notation for general functions $g_u$ with a non-increasing upper tail with threshold $d_u$ ends the proof.

\subsubsection*{Expressions \eqref{T-h-31}--\eqref{T-h-42}}
The proof is similar to that of \eqref{T-h-21} and \eqref{T-h-22}, where we can use Lemma \ref{Lemma-Link}.

\subsection{Proof of Lemma \ref{LemmaLink2}}\label{Proof:LemmaLink2}
We prove the result for
$$(X_1,X_2) \text{ upper-upper comonotonic } \Longleftrightarrow (-X_1,-X_2) \text{ lower-lower comonotonic},$$
and similar arguments can be used to prove the other results.\\

Let $\underline{X}=\left(X_1,X_2\right)$ be an upper-upper comonotonic random vector with threshold $\left(a_1,a_2\right)$, and we denote its support by $\mathcal{S}$. We consider the random vector $\underline{\bar{X}}=\left(-X_1,-X_2\right)$ with support $\bar{\mathcal{S}}$. For any $(x_1,x_2)\in \mathcal{S}$, we have that $\left(-x_1,-x_2\right) \in \bar{\mathcal{S}}$, which leads to:
$$
\left(x_1,x_2\right) \in \mathcal{S} \cap TR(\left(a_1,a_2\right)) \Longleftrightarrow\left(-x_1,-x_2\right)\in \bar{\mathcal{S}}\cap BL((-a_1,-a_2)),
$$	
and therefore:
$$
\mathcal{S}\cap TR((a_1,a_2))\ \text{is a comonotonic set}\Longleftrightarrow\bar{\mathcal{S}}\cap BL((-a_1,-a_2))\ \text{is a comonotonic set}.
$$
The proof of the equivalence between the second conditions for lower-lower and upper-upper comonotonicity follows directly as:
$$\mathbb{P}\left[X_1<a_1,X_2<a_2\right]>0  \Leftrightarrow \mathbb{P}\left[\bar{X}_1>-a_1,\bar{X}_2>-a_2\right]>0.$$

For the third condition, it is also straightforward to show that:
$$\bar{\mathcal{S}} \cap \left(\mathbb{R}^2\setminus\left(TR(\underline{a})\cup\overline{BL}(\underline{a})\right)\right) \text{\ is empty\ }\Longleftrightarrow \bar{\mathcal{S}}\cap \left(\mathbb{R}^2\setminus\left(BL(-\underline{a})\cup\overline{TR}(-\underline{a})\right)\right) \text{\ is empty\ } $$

Combining the equivalences between the three conditions ends the proof.

\subsection{Proof of Theorem \ref{L-1}}\label{Proof:L-1}
\subsubsection*{Fisrt statement: proof of $(\Rightarrow)$}
Suppose that $\left(X_2,X_1-X_2\right)$ is upper-upper comonotonic with threshold $\underline{a}$ given by \eqref{UpperComonotonic} with $\alpha\in[0,1]$, and consider the function $h_\alpha$ defined in \eqref{H_Alpha}. We want to show that the function $h_\alpha$ has a non-decreasing upper tail with threshold $\pi$. From the definition of non-decreasing upper tail functions, we need to show that $h_\alpha(p)\leq h_\alpha(\pi)$ for $p<\pi$, and that $h_\alpha(p_1)\leq h_\alpha(p_2)$, for $\pi \leq p_1 \leq p_2$.

Notice that because $\left(X_1,X_2\right)$ is comonotonic, we have that $\left(X_1,X_2\right) \overset{d}{=}\left(F^{-1(\alpha)}_{X_1}(U),F^{-1(\alpha)}_{X_2}(U)\right)$, where $U$ is uniformly distributed over $[0,1]$. Hence, $\mathbb{P}\left[\left(X_2,X_1-X_2\right)\in \mathcal{S}\right]=1$, with:
\begin{equation}\label{L-1_5}\mathcal{S}=\left\{\left(F^{-1(\alpha)}_{X_2}(p),h_\alpha(p)\right)|p\in[0,1]\right\}.\end{equation}
Moreover, since $\left(X_2,X_1-X_2\right)$ is upper-upper comonotonic, then $\mathcal{S}\cap \left(\mathbb{R}^2\setminus\left(TR(\underline{a})\cup\overline{BL}(\underline{a})\right)\right)$ is empty, i.e. the support $\mathcal{S}$ takes values only in the top-right quadrant $TR(\underline{a})$ and the bottom-left quadrant $\overline{BL}(\underline{a})$.

We now prove the first property, i.e. $h_\alpha(p)\leq h_\alpha(\pi)$ for $p<\pi$. Let $p< \pi$ and consider the element $\left(F^{-1(\alpha)}_{X_2}(p),h_\alpha(p)\right)\in \mathcal{S}$. Since quantiles are non-decreasing, we have that $F^{-1(\alpha)}_{X_2}(p)\leq F^{-1(\alpha)}_{X_2}(\pi)$. Moreover, $\mathcal{S}\cap \left(\mathbb{R}^2\setminus\left(TR(\underline{a})\cup\overline{BL}(\underline{a})\right)\right)$ is empty, which means that the element $\left(F^{-1(\alpha)}_{X_2}(p),h_\alpha(p)\right)$ is necessarily in the bottom-left quadrant, and hence:
\begin{equation}\label{L-1_3} h_\alpha(p)\leq h_\alpha(\pi), \qquad \text{for } p< \pi.\end{equation}

We now prove the second property, i.e. $h_\alpha(p_1)\leq h_\alpha(p_2)$, for $\pi \leq p_1 \leq p_2$. Consider the elements $p_1$ and $p_2$ with $\pi \leq p_1 \leq p_2$. Because $F^{-1(\alpha)}_{X_2}$ is non-decreasing, we have that $F^{-1(\alpha)}_{X_2}(p_j)\geq F^{-1(\alpha)}_{X_2}(\pi)$ for $j=1,2$. Note that the top-right quadrant $TR(\underline{a})$ does not contain the element $\underline{a}$. Thus, we distinguish two cases. In the first case where $F^{-1(\alpha)}_{X_2}(p_j)=F^{-1(\alpha)}_{X_2}(\pi)$, we have that $h_\alpha(p_j)=F^{-1(\alpha)}_{X_1}(p_j)-F^{-1(\alpha)}_{X_2}(\pi)$, and by the monotonicity of $F^{-1(\alpha)}_{X_1}$, we find that $h_\alpha(p_j)\geq h_\alpha(\pi)$. In the second case where $F^{-1(\alpha)}_{X_2}(p_j)>F^{-1(\alpha)}_{X_2}(\pi)$, we use again the fact that $\left(F^{-1(\alpha)}_{X_2}(p_j),h_{\alpha}(p_j)\right)$ is necessarily in the top-right quadrant, and hence, $h_\alpha(p_j)> h_\alpha(\pi)$ for $j=1,2$ also holds. Thus, we have that $h_\alpha(p_j)\geq h_\alpha(\pi)$ for $j=1,2$ and $\pi\leq p_1\leq p_2$.

Invoking again the non-decreasingness of $F^{-1(\alpha)}_{X_2}$, together with the monotonicity of the set $\mathcal{S}_{\underline{a}} = \mathcal{S}\cap TR(\underline{a})$, leads to:
\begin{equation}\label{L-1_4}h_\alpha(p_1)\leq h_\alpha(p_2), \qquad \text{for } \pi\leq p_1 \leq p_2.\end{equation}
Inequalities \eqref{L-1_3} and \eqref{L-1_4} imply that $h_\alpha$ has a non-decreasing upper tail with threshold $\pi$.
\subsubsection*{First statement: proof of $(\Leftarrow)$}
Suppose the function $h_\alpha$ defined in \eqref{H_Alpha} has a non-decreasing upper tail with threshold $\pi\in(0,1)$. We aim at showing that the vector $(X_2,X_1-X_2)$ is upper-upper comonotonic with threshold $\underline{a}$ given in \eqref{UpperComonotonic}. For this, we need to show that the three conditions for upper-upper comonotonicity given in Definition \ref{Def:UpperCom} are satisfied. Namely, that $\mathcal{S}_{\underline{a}}$ is a comonotonic set, that $\mathbb{P}[\underline{X}>\underline{a}]>0$, and that the set $\mathcal{S}\cap \left(\mathbb{R}^2\setminus\left(TR(\underline{a})\cup\overline{BL}(\underline{a})\right)\right)$ is empty.

Recall that the support of $(X_2,X_1-X_2)$ is given by the set $\mathcal{S}$ defined in \eqref{L-1_5}. Moreover, the point $\underline{a}=\left(F^{-1(\alpha)}_{X_2}(\pi),h_\alpha(\pi)\right)$ is an element of $\mathcal{S}$.

First, the set $\mathcal{S}_{\underline{a}}$ is given by:
$$\mathcal{S}_{\underline{a}} = \mathcal{S}\cap TR(\underline{a}) = \left\{\left(F^{-1(\alpha)}_{X_2}(p),h_\alpha(p) \right)| p>\pi\right\}.$$
The quantile $F^{-1(\alpha)}_{X_2}$ is non-decreasing for $p\in[0,1]$, and $h$ has a non-decreasing upper tail with threshold $\pi$, and hence, $h_\alpha$ is non-decreasing for $p\geq \pi$. Thus, $\mathcal{S}_{\underline{a}}$ is a comonotonic set.

Second, for $\pi<1$, the set $\mathcal{S}_{\underline{a}}$ is not empty, and hence $\mathbb{P}[\left(X_2,X_1-X_2\right)> \underline{a}]>0$.

Third, we combine again the fact that the quantile $F^{-1(\alpha)}_{X_2}$ is non-decreasing and that $h_\alpha$ has a non-decreasing upper tail with threshold $\pi$ to conclude that, on the one hand, we have $h_\alpha(p)\geq h_\alpha(\pi)$ and $F^{-1(\alpha)}_{X_2}(p)\geq F^{-1(\alpha)}_{X_2}(\pi)$ for $p\geq \pi$, and on the other hand, we have $h_\alpha(p)\leq h_\alpha(\pi)$ and $F^{-1(\alpha)}_{X_2}(p)\leq F^{-1(\alpha)}_{X_2}(\pi)$ for $p< \pi$. Therefore, the support $\mathcal{S}$ of $\left(X_2,X_1-X_2\right)$ takes values in either the top-right or the bottom-left quadrants. We conclude that the set $\mathcal{S}\cap\left(\mathbb{R}^2\setminus\left(TR(\underline{a})\cup\overline{BL}(\underline{a})\right)\right)$ is empty, which proves the result.

\subsubsection*{Second statement}
For $\pi\in (0,1)$ and $\alpha \in [0,1]$, let $\left(X_2,X_1-X_2\right)$ be lower-lower comonotonic with threshold $$\underline{a}=\left(a_1,a_2\right)=\left(F^{-1(\alpha)}_{X_2}(\pi),F^{-1(\alpha)}_{X_1}(\pi)-F^{-1(\alpha)}_{X_2}(\pi)\right),$$
Lemma \ref{LemmaLink2} states that the random vector $\left(\bar{R}_2,\bar{R}_1-\bar{R}_2\right)$, where $\bar{R}_i=-X_i$, is upper-upper comonotonic with threshold $$\underline{\bar{a}}=\left(-a_1,-a_2\right)=\left(F^{-1(1-\alpha)}_{\bar{R}_2}(1-\pi),F^{-1(1-\alpha)}_{\bar{R}_1}(1-\pi)-F^{-1(1-\alpha)}_{\bar{R}_2}(1-\pi)\right),$$
where we use the fact that $-F^{-1(\alpha)}_{X_i}(\pi)=F^{-1(1-\alpha)}_{\bar{R}_i}(1-\pi)$. 
We define the function $\bar{h}_{\alpha}$ as follows:
$$\bar{h}_{\alpha}(p) = F^{-1(\alpha)}_{\bar{R}_1}(p)-F^{-1(\alpha)}_{\bar{R}_2}(p),\qquad \text{for }p\in(0,1).$$
We can write:
$$\underline{\bar{a}}=\left(F^{-1(1-\alpha)}_{\bar{R}_2}(1-\pi),\bar{h}_{1-\alpha}(1-\pi)\right).$$
Since $\left(X_1,X_2\right)$ is comonotonic, we also have that $\left(\bar{R}_1,\bar{R}_2\right)$ is comonotonic. Thus, from Lemma \ref{L-1}, the upper-upper comonotonicity of $\left(X_2,X_1-X_2\right)$ is equivalent with $\bar{h}_{1-\alpha}$ having a non-decreasing upper tail with threshold $1-\pi$.

Define the function $h_\alpha$ as follows:
$$h_\alpha(p)=F^{-1(\alpha)}_{X_1}(p) - F^{-1(\alpha)}_{X_2}(p).$$
We have that $h_\alpha(p)=-\bar{h}_{1-\alpha}(1-p)$ for $p\in (0,1)$. Using Lemma \ref{Lemma-Link}, the function $-\bar{h}_{1-\alpha}(-p)$, for $p\in(-1,0)$, has a non-decreasing lower tail with threshold $\pi-1$. A change of variable leads to the conclusion that $h(p)=-\bar{h}_{1-\alpha}(1-p)$, for $p\in(0,1)$, has a non-decreasing lower tail with threshold $\pi$, which proves the result.

\subsection{Proof of Proposition \ref{LN-case2-functiong}}\label{LN-case2-proof}
	We start with the monotonicity of the function $h$. The cdf $F_{W}$ is strictly increasing and continuous, and hence, its inverse $F_W^{-1}$ is differentiable. This implies that the function $h$ defined in \eqref{LNg} is also differentiable. We have that $$\frac{\text{d}F_W^{-1}(u)}{\text{d}u}=\frac{1}{f_W\left(F_W^{-1}(u) \right)}>0,$$
which leads to the following expression for the derivative of $h$:
$$\frac{\text{d}h}{\text{d}u}=\frac{1}{f_W\left(F_W^{-1}(u) \right)}\left(\sigma_1\text{e}^{\mu_1+\sigma_1F_W^{-1}(u)}-\sigma_2\text{e}^{\mu_2+\sigma_2F_W^{-1}(u)} \right).$$
Studying the sign of this derivative is equivalent to studying the sign of
$$\log \frac{\sigma_1}{\sigma_2}+ \mu_1-\mu_2 - (\sigma_2-\sigma_1)F_W^{-1}(u).$$
If $\sigma_1=\sigma_2$, the function $h$ is always monotone, and the sign of $\frac{\text{d}h}{\text{d}u}$ depends on the sign of $\mu_1-\mu_2$. If $\sigma_1<\sigma_2$, we find that $\frac{\text{d}h}{\text{d}u}$ is positive on $[0,u^{\star}]$ and negative on $[u^{\star},1]$, where 
\begin{equation}\label{ustar}
u^{\star}=
F_W\left( \frac{\log \frac{\sigma_1}{\sigma_2}+\mu_1-\mu_2}{\sigma_2-\sigma_1}\right).
\end{equation}
Hence, the function $h$ starts at zero, is non-decreasing until it reaches its maximum in $ u^{\star}$, after which the function $h$ is non-increasing and tends to $-\infty$. Analogously, if $\sigma_1>\sigma_2$, we find that the function $h$ starts at zero and is non-increasing in $[0,u^{\star}]$. In the interval $[u^{\star},1]$, the function is non-decreasing and tends to $+\infty$. Also, we find that for $\sigma_1 \neq \sigma_2$, the function $h$ has a monotone upper tail with threshold $c^{\star}$ such that $h(c^{\star})=0$, and hence, $c^{\star}=F_W\left(\frac{\mu_2-\mu_1}{\sigma_1-\sigma_2} \right)$.

\end{document}